\documentclass[aap,MSNbibl,dvips]{arximspdf}
\usepackage{graphicx}
\doi{10.1214/13-AAP944} %kopijuoti is PTS
\volume{24}
\issue{3}
\pubyear{2014}
\firstpage{1172}
\lastpage{1198}

\makeatletter
\newcommand{\rrvert}{\vert}
\newcommand{\llvert}{\vert}
\newtheorem{theorem}{Theorem}[section]
\newtheorem{lemma}[theorem]{Lemma}
\newtheorem{proposition}[theorem]{Proposition}
\newtheorem{corollary}[theorem]{Corollary}
\newproclaim{remark}[theorem]{Remark}
\newproclaim{definition}[theorem]{Definition}
\makeatother

\begin{document}
\begin{frontmatter}

\title{Moments and Lyapunov exponents for\\ the parabolic Anderson model}
\runtitle{Moments and Lyapunov exponents}

\begin{aug}
\author[A]{\fnms{Alexei} \snm{Borodin}\ead[label=e1]{borodin@math.mit.edu}\thanksref{t1}}
\and
\author[B]{\fnms{Ivan} \snm{Corwin}\corref{}\ead[label=e2]{icorwin@mit.edu}\thanksref{t2}}
\runauthor{A. Borodin and I. Corwin}
\affiliation{Massachusetts Institute of Technology and Institute for
Information Transmission Problems, and
Massachusetts Institute of Technology and\\ Clay Mathematics Institute}
\address[A]{Department of Mathematics\\
Massachusetts Institute of Technology\\
77 Massachusetts Avenue\\
Cambridge, Massachusetts 02139-4307\\
USA\\
and\\
Institute for Information\\
\quad Transmission Problems\\
Bolshoy Karetny per. 19\\
Moscow 127994\\
Russia\\
\printead{e1}} %adresu isvedimo komanda gale!
\address[B]{Department of Mathematics\\
Massachusetts Institute of Technology\\
77 Massachusetts Avenue\\
Cambridge, Massachusetts 02139-4307\\
USA\\
and\\
Clay Mathematics Institute\\
10 Memorial Blvd. Suite 902\\
Providence, Rhode Island 02903\\
USA\\
\printead{e2}}
\end{aug}
\thankstext{t1}{Supported in part by NSF Grant DMS-10-56390.}
\thankstext{t2}{Supported in part by the NSF through PIRE Grant OISE-07-30136 and DMS-12-08998
as well as by Microsoft Research through the Schramm Memorial
Fellowship, and by the Clay Mathematics Institute through a Clay
Research Fellowship.}

% HISTORY:
\received{\smonth{1} \syear{2013}}
\revised{\smonth{6} \syear{2013}}

% ABSTRACT
%
\begin{abstract}
We study the parabolic Anderson model in $(1+1)$ dimensions with
nearest neighbor jumps and space--time white noise (discrete
space/continuous time). We prove a contour integral formula for the
second moment and compute the second moment Lyapunov exponent. For the
model with only jumps to the right, we prove a contour integral formula
for all moments and compute moment Lyapunov exponents of all orders.
\end{abstract}

% KEYWORDS

\end{frontmatter}

%s1 #&#
\section{Introduction and main results}\label{sec1}
%s1.1 #&#
\subsection{Nearest-neighbor parabolic Anderson model}
The nearest-neighbor pa\-ra\-bolic Anderson model on $\mathbb{Z}$ is
the solution to a coupled system of diffusions on $[0,\infty)$ given by
%
%
%e1 #&#
\begin{equation}
\label{notmild} \frac{d}{dt} Z_{\beta}(t,n) = \frac{1}{2}
\Delta^{p,q}Z_{\beta
}(t,n) + \beta Z_{\beta}(t,n)
\,dW_n(t).
\end{equation}
We focus here on delta function initial data $Z_{\beta}(0,n) = {\mathbf
1}_{n=0}$. Here $t\in\mathbb{R}_{+}$, $n\in\mathbb{Z}$, and the
operator $\Delta ^{p,q}$ (which is the generator for a nearest neighbor
continuous time random walk) acts on functions $f(n)$ as
%
%
%e2 #&#
\begin{equation}
\label{deltadef} \Delta^{p,q} f (n) = pf(n-1) + qf(n+1) - (p+q) f(n).
\end{equation}
We assume that $p,q\geq0$ and $p+q=2$. The collection $\{W_n(\cdot)\}
_{n\in\mathbb{Z}}$ are independent Brownian motions and $\beta\in
\mathbb{R}_{+}$.
%%For a general physical and mathematical background on this model, the
%reader is referred to \cite{CarMol,denHoll} as well as part I of

%s1.1.1 #&#
\subsubsection{Population growth in random environment}
The coupled diffusions can be considered as modeling population growth
in a random, quickly changing environment at each spatial location, and
with migration between locations. Consider a population of many small
particles living on the sites of $\mathbb{Z}$. There are three forces acting
upon this system:
\begin{longlist}[(3)]
\item[(1)] Each particle at time $t$ and lattice site $n$ independently
    duplicates itself at rate $r_{+}(t,n)$;

\item[(2)] Each particle at time $t$ and lattice site $n$ independently
    dies at rate $r_{-}(t,n)$;

\item[(3)] Each particle at time $t$ and lattice site $n$ independently
    jumps to a neighboring site $n-1$ with rate $q/2$ and $n+1$ with
    rate $p/2$.
\end{longlist}
Letting $m(t,n)$ be the expected population size at time $t$ and
location $n$, one finds that \cite{CarMol}
\[
\frac{d}{dt} m(t,n) = \frac{1}{2}\Delta^{p,q} m(t,n) +
\bigl(r_{+}(t,n)-r_{-}(t,n) \bigr) m(t,n).
\]
If the duplication and death rates are independent in space and quickly
mixing in time, the factor $(r_{+}(t,n)-r_{-}(t,n))$ is well modeled by
$\beta \,dW_n(t)$ where $\beta$~modulates the relative rates of jumping
and duplication/death. The delta function initial data translates into
starting with all the particles clustered at the origin and then
allowing them to spread over time.

As explained in \cite{CarMol,denHoll}, it is of physical interest for
these models to understand the structure of regions in space--time in
which the population size is significantly larger than expected. This
phenomenon is called \textit{intermittency}. Generally, one seeks to
measure the effect of changing various parameters with respect to this
phenomenon. Of specific interest are the spatial dimension (replacing
$\mathbb{Z}$ by~$\mathbb{Z}^d$), the strength of $\beta$, and the type
of environmental noise (replacing space--time noise by spatially
varying but constant in time noise, or noise which is itself built out
of interacting particle systems); see part~I of \cite{Fest} and
\cite{CC} for reviews of these various directions. In the present paper
we restrict ourself to the one-dimensional, space--time independent
case and offer a new approach to computing the moments of this model.
Section~\ref{sec12} below explains the relevance of the moments to the
intermittency phenomenon.

%s1.1.2 #&#
\subsubsection{Directed polymers}\label{dps}
Closely related to the above branching diffusion representation, the
Feynman--Kac representation for this coupled system of diffusions
writes $Z_{\beta}(t,n)$ as point to point partition functions for a
random polymer model
%
%
%e3 #&#
\begin{equation}
\label{FK} Z_\beta(t,n) = \mathcal{E}_{\pi(0)=0} \biggl[
\mathbf{1}_{\pi
(t)=n} \exp \biggl( \int_0^t
\beta \,dW_{\pi(s)}(s) \,ds - \frac{\beta
^2 t}{2} \biggr) \biggr],
\end{equation}
where $\pi(s)$ is a Markov process with state space $\mathbb{Z}$ and
generator given by $\frac{1}{2}\Delta^{q,p}$ (which is the adjoint of
$\frac{1}{2}\Delta^{p,q}$), and $\mathcal{E}_{\pi(0)=0}$ is the
expectation with respect to starting $\pi(0)=0$. We write $\mathbb {E}$
for the expectation over the disorder. The polymer measure on paths
$\pi (\cdot)$ is defined as the argument of the above expectation,
normalized by $Z_{\beta}(t,n)$.

Directed polymer models are important from a number of perspectives;
see \mbox{\cite{CSY,ICReview}} and references therein. They were
introduced to study the domain walls of Ising models with impurities at
high temperature and have been applied to other problems like vortices
in superconductors, roughness of crack interfaces, Burgers turbulence,
and interfaces in competing bacterial colonies. They also provide a
unified mathematical framework for studying a variety of different
abstract and physical problems including some in stochastic
optimization, bio-statistics, queuing theory and operations research,
interacting particle systems and random growth models.

The above defined class of directed polymers is a generalization of the
model (at $p=2$ and $q=0$) introduced by O'Connell--Yor
\cite{OCon-Yor}. The primary interest in the study of directed polymers
is to understand the free energy fluctuations [i.e., $\log Z_{\beta
}(t,n)$] and the transversal path fluctuations of the polymer measure
under the limit at $t$ and $n$ go to infinity. For the special $p=2$
and $q=0$ case, there has been significant progress in both of these
directions coming from the work of
\mbox{\cite{OCon,SeppValko,BorCor,BCF}}. The model is now known to be
in the Kardar--Parisi--Zhang universality class, which predicts these
asymptotic fluctuation behaviors. It is expected that this asymptotic
behavior should not depend on the values of $p,q$ and $\beta$. The
present analysis of the moments of $Z_{\beta }(t,n)$ constitute a step
toward an analysis of this class. On the other hand, one should note
that by virtue of the intermittency which we prove herein, one knows
that these moments will not determine the distribution of
$Z_{\beta}(t,n)$.

From the above polymer representation for $Z_{\beta}(t,n)$ and the
Gaussian nature of the noise, one sees (by interchanging the path
expectations with the expectation over the disorder) that
\begin{eqnarray*}
&& \mathbb{E} \Biggl[\prod_{i=1}^2 Z_{\beta}(t,n_i) \Biggr]
\\
&&\qquad =
\mathcal{E}_{\pi _1(0)=\pi_2(0)=0} \biggl[
\mathbf{1}_{\pi_1(t)=n_1,\pi_2(t)=n_2} \exp \biggl(\frac{\beta ^2}{2}
\int_0^t \mathbf{1}_{\pi_1(s)=\pi _2(s)} \,ds \biggr) \biggr].
\end{eqnarray*}
In other words, letting $\pi= \pi_1-\pi_2$ and $\mathcal{E}$ be the
associated expectation, we find that
\[
\mathbb{E} \Biggl[\prod_{i=1}^2
Z_{\beta}(t,n_i) \Biggr] = \mathcal{E} \biggl[
\mathbf{1}_{\pi(t)=n_1-n_2} \exp \biggl(\frac{\beta^2}{2} \int_0^t
\mathbf{1}_{\pi(s)=0} \,ds \biggr) \biggr].
\]
This is the first moment of the partition function for a random walk
$\pi$ which feels a pinning potential of strength $\beta^2/2$ at the
origin; see \cite{Berger} and references therein for more discussion on
this model.

%s1.2 #&#
\subsection{Lyapunov exponents and intermittency}\label{sec12}
In order to introduce and explain the mathematical definition of
intermittency, we introduce two types of the Lyapunov exponents for the
parabolic Anderson model. Consider a velocity $\nu\in\mathbb{R}$. Then
the \textit{almost sure Lyapunov exponent} with respect to velocity
$\nu$ is given by
%
%
%e4 #&#
\begin{equation}
\label{tildegamma} \tilde\gamma_1(\beta;\nu) = \lim_{t\to\infty}
\frac{1}{t} \log Z_{\beta} \bigl(t,\lfloor\nu t\rfloor \bigr).
\end{equation}
The existence of this almost sure limit is due to a sub-additivity
argument (see \cite{CarMol}, Section~IV.1). The \textit{$p$th moment
Lyapunov exponent} with respect to velocity $\nu$ is given by
%
%
%e5 #&#
\begin{equation}
\label{gamma} \gamma_k(\beta;\nu) = \lim_{t\to\infty}
\frac{1}{t} \log\mathbb{E} \bigl[ \bigl(Z_{\beta} \bigl(t,\lfloor\nu
t \rfloor \bigr) \bigr)^k \bigr].
\end{equation}

If the initial data $Z_{\beta}(0,n)$ is stationary with respect to
shifts in $n$, then the exponents are, in fact, independent of $\nu$.
We, however, consider initial data in which $Z_{\beta}(0,n) = \mathbf
{1}_{n=0}$, and hence the exponents will depend on the velocity $\nu$
nontrivially.
%
%
%de1.1 #&#
\begin{definition}
A parabolic Anderson model shows intermittency if the Lypanov exponents
are strictly ordered as
\[
\tilde\gamma_1(\beta;\nu) < \gamma_1(\beta;\nu) <
\frac{\gamma
_{2}(\beta;\nu)}{2} <\frac{\gamma_{3}(\beta;\nu)}{3} < \cdots.
\]
\end{definition}
The weak ordering of exponents is a consequence of Jensen's inequality
(for the first inequality) and H\"{o}lder's inequality (for all
subsequent inequalities). A useful fact is recorded in the following
(cf. \cite{CarMol}, Theorem III.1.2):

%
%le1.2 #&#
\begin{lemma}\label{inter}
If for any $k\geq1$,
%
%
%e6 #&#
\begin{equation}
\label{base1} \frac{\gamma_k(\beta;\nu)}{k}<\frac{\gamma_{k+1}(\beta;\nu)}{k+1}
\end{equation}
then for all $p\geq k$
%
%
%e7 #&#
\begin{equation}
\label{base2} \frac{\gamma_p(\beta;\nu)}{p}<\frac{\gamma_{p+1}(\beta;\nu)}{p+1}.
\end{equation}
\end{lemma}
%
%H\"{o}lder's inequality implies that
%$$\EE\left[(Z_{\beta}(t,\lfloor\var t\rfloor))^k\right]^2 \leq\EE
%from which follows the inequality
%If (\ref{base1}) holds, then
%$$\gamma_k(\beta;\var) < \frac{k}{k+1} \gamma_{k+1}(\beta;\var).$$
%Plug this into the right-hand side (\ref{hineq}) with $h=1$ and $k$
%replaced by $k+1$. Rearranging terms yields (\ref{base2}) for $p=k+1$.
%This can be repeated inductively, and yields the claimed result.

As explained in \cite{CarMol}, intermittent random fields are
distinguished by the formation of strong pronounced spatial structures
such as sharp peaks which give the main contribution to the physical
processes in such media. A popular example cited therein is the
observation of Zeldovich that the Solar magnetic field is intermittent
since more than 99$\%$ of the magnetic energy concentrates on less than
1$\%$ of the surface area.\vadjust{\goodbreak}

The above mathematical definition of intermittency is related to the
presence of high peaks of $Z_{\beta}(t,n)$ that dominate large time
moment asymptotics. In particular, fix $\alpha$ such that $\frac
{\gamma_{k}(\beta;\nu)}{k} < \alpha<
\frac{\gamma_{k+1}(\beta;\nu)}{k+1}$; then we know that
$\mathbb{P}(Z_{\beta}(t,\nu t) > e^{\alpha t}) >0$ as $t\to\infty$.
Writing
\begin{eqnarray*}
&& \mathbb{E} \bigl[ Z_{\beta} \bigl(t,\lfloor\nu t \rfloor
\bigr)^{k+1} \bigr]
\\
&&\qquad  = \mathbb{E} \bigl[ Z_{\beta} \bigl(t,\lfloor \nu t
\rfloor \bigr)^{k+1}\mathbf{1}_{Z_{\beta }(t,\lfloor\nu
t\rfloor)<e^{\alpha t}} \bigr] + \mathbb{E} \bigl[ Z_{\beta }
\bigl(t,\lfloor\nu t\rfloor \bigr)^{k+1}
\mathbf{1}_{Z_{\beta}(t,\lfloor\nu t\rfloor)>e^{\alpha t}} \bigr],
\end{eqnarray*}
we observe that the first term is \mbox{${\leq}e^{\alpha(k+1) t}$}, but
the sum of the two is asymptotically $e^{\gamma_{k+1}(\beta;\nu) t}$,
which is exponentially (as $t$ grows) larger than $e^{\alpha(k+1)t}$.
This means that the event $ \{Z_{\beta}(t,\nu t) > e^{\alpha t} \}$
gives overwhelming contribution to the $(k+1)$st moment. On the other
hand,
\[
\mathbb{E} \bigl[ Z_{\beta} \bigl(t,\lfloor\nu t\rfloor
\bigr)^k \bigr] \geq e^{\alpha k t} \mathbb{P} \bigl(
Z_{\beta} \bigl(t,\lfloor\nu t \rfloor \bigr) > e^{\alpha t} \bigr)
\]
and hence, for large $t$
\[
\mathbb{P} \bigl( Z_{\beta} \bigl(t,\lfloor\nu t\rfloor \bigr) >
e^{\alpha
t} \bigr) \leq\frac{e^{\gamma_k(\beta;\nu) t}}{e^{\alpha k t}} =\exp \biggl\{ - \biggl(\alpha-
\frac{\gamma_k(\beta;\nu)}{k} \biggr) t \biggr\},
\]
which is exponentially small.

In the case of spatially translation invariant ergodic solutions
$Z_{\beta}(t,n)$, the consequences of intermittency may be interpreted
via spatial averages over large balls at a fixed (large) time. Thus one
can talk about islands where the solution is at least $e^{((\gamma
_k(\beta;\nu))/k) t}$ (as opposed to the typical value of $e^{\tilde
\gamma_1(\beta;\nu)t}$) whose spatial density is not more than $e^{-
(((\gamma_{k+1}(\beta;\nu))/(k+1)) -((\gamma _k(\beta;\nu))/k)) t}$.
Our results are for delta initial data $Z_{\beta}(0,n) =
\mathbf{1}_{n=0}$ and not stationary initial data.

In terms of the population model interpretation of $Z_{\beta}(t,n)$,
the above discussion implies that knowledge of the Lyaponov exponents
translates into detailed information about the spatial frequency of
large clusters of population growth in space.

In this direction, the primary contribution of this paper is the
precise calculation of the first two moment Lyapunov exponents for the
general $p,q$ model and the calculation of all moment Lyapunov
exponents for the special $p=2$ and $q=0$ case. From the above
considerations, this provides detailed information about the
intermittent structure of the corresponding population growth models.

%The inequality between $\tilde\gamma_1$ and $\gamma_1$ implies that
%the typical growth rate of $Z_{\beta}(t,\var t)$ is small than the
%growth rate of its first moment, thus indicating that with probability
%going to zero (as $t\to\infty$), there exist sufficiently high values
%of $Z_{\beta}(t,\var t)$

%s1.3 #&#
\subsection{Main results}

All of our results pertain to the nearest-neighbor para\-bolic Anderson
model with delta initial data: $Z_{\beta}(0,n) = \mathbf{1}_{n=0}$.
Our first result is a formula for the two-point moment of the
model.

%
%th1.3 #&#
\begin{theorem}\label{cor1}
For $n_1\geq n_2$,
%
%
%e8 #&#
\begin{eqnarray}\label{2momentformula}
\mathbb{E} \Biggl[\prod_{i=1}^{2}
Z_\beta(t,n_i) \Biggr] &=& \frac
{1}{(2\pi\iota)^2} \oint\oint
\frac{ (pz_1-qz_1^{-1}) - (pz_2 - q z_2^{-1})}{
(pz_1-qz_1^{-1}) - (pz_2 - q z_2^{-1}) - 2\beta^2}
\nonumber\\[-8pt]\\[-8pt]
&&\hspace*{54pt}{}\times  F^{p,q}_{t,n_1}(z_1)
F^{p,q}_{t,n_2}(z_2) \frac{dz_1}{z_1}
\frac{dz_2}{z_2},\nonumber
\end{eqnarray}
where
\[
F^{p,q}_{t,n}(z) = z^{-n} e^{(t/2)(pz+qz^{-1} -2)}
\]
and where the contour of $z_1$ is the unit circle, and the contour for
$z_2$ is a circle around 0 of radius sufficiently small so as not to
include any poles of the integrand aside from $z_2=0$.
\end{theorem}
This theorem is proved in Section~\ref{nnpamsec}. This also provides an
exact formula for the first moment of the pinned polymer partition
function discussed above in Section~\ref{dps}. The following corollary
follows immediately from the exact result above along with Lemma
\ref{inter} applied for $k=1$. In the symmetric case $p=q$ this was
established as a special case of the results in \cite{CarMol},
Chapter~III.
%
%
%co1.4 #&#
\begin{corollary}
The $p,q$ nearest-neighbor parabolic Anderson model displays
intermittency at the velocity $p-q$.
\end{corollary}
%
%We state an immediate corollary of this theorem which shows how it
%yields an exact solution to a renewal equation. Let $p(t,n)$ denote
%the heat kernel for the generator $\Delta^{p,q}$ (i.e., the solution
%to $\frac{d}{dt} p(t,n) = \Delta^{p,q} p(t,n)$ with $p(0,n) =
%For all $t\in\Rplus$, $n\in\Z$, and $\beta\in\Rplus$ the equation
%$$f(t,n)^2 = p(t,n)^2 + \beta^2 \int_{0}^{t} \,ds \sum_{m=-\infty}^{
%with initial data $f(0,n)=\mathbf{1}_{n=0}$ is solved by the
%right-hand side of (\ref{2momentformula}) with $n_1=n_2=n$.
%The integral formulation of (\ref{notmild}) is
%$$
%Z_{\beta}(t,n) = p(t,n) + \int_0^t \sum_{m=-\infty}^{\infty}
%p(t-s,n-m) Z_{\beta}(s,m) \,dW_m(s).
%$$
%Square both sides and take expectations. Call $f(t,n) = \EE\left[Z_{
%result.

Via asymptotic analysis, Theorem \ref{cor1} enables us to calculate
the first and second moment Lyapunov exponents for the parabolic
Anderson model (as well as the first moment Lyapunov for the pinned
polymer partition function).

%
%th1.5 #&#
\begin{theorem}\label{secondlyap}
The first moment Lyapunov exponent at velocity $p-q$, for the
nearest-neighbor parabolic Anderson model is given by $\gamma_1(\beta;p-q)=0$.
The second moment Lyapunov exponent at velocity $p-q$, for the
nearest-neighbor parabolic Anderson model is given by
\[
\gamma_2(\beta;p-q) = H_2 \bigl(z^0_2
\bigr),
\]
where
\begin{eqnarray*}
H_2(z) &=&\tfrac{1}{2} \bigl( p s(z)+q \bigl(s(z) \bigr)^{-1} -2
-(p-q)\log \bigl(s(z) \bigr)
\\
&&\hspace*{51pt}{} + p z +q z^{-1} -2 -(p-q)\log z \bigr)
\end{eqnarray*}
with
\[
s(z) = \frac{(pz-qz^{-1}+2\beta^2) + \sqrt{(pz-qz^{-1}+2\beta^2)^2
+ 4pq}}{2p}
\]
and where $z^0_2$ is the unique solution to $H_2'(z)=0$ over $z\in
(0,\infty)$.\vadjust{\goodbreak}

When $p=q=1$,
\begin{eqnarray*}
z^0_2 &=& \tfrac{1}{2} \bigl(-\beta^2+ \sqrt{4+\beta^4} \bigr),
\qquad s(z_0)= \tfrac{1}{2} \bigl(\beta^2+ \sqrt{4+\beta^4} \bigr),
\\
H_2 \bigl(z^0_k \bigr) &=& 2 \bigl(\sqrt{4+ \beta^4}-2 \bigr),
\end{eqnarray*}
which implies that
\[
\gamma_2(\beta;0) = 2 \bigl(\sqrt{4+\beta^4}-2 \bigr)
\]
for the standard ($p=q$) parabolic Anderson model.
\end{theorem}
This theorem is proved in Section~\ref{nnpamsec} via asymptotic
analysis of Theorem \ref{cor1}. We include the full details only for
the case $p=q=1$.

%
%re1.6 #&#
\begin{remark}
The above theorem is stated only for a velocity given by $\nu=p-q$.
For general $p-q\neq0$ the same approach as given in the proof
provides the exact values of the first and second moment Lyapunov
exponents, but we forgo including this herein.
\end{remark}

%attainable via studying the local time at zero of the random walk
%performed by the difference of two paths}

We now turn our attention to the one-sided case of the nearest-neighbor
parabolic Anderson model, where $p=2$ and $q=0$. In this case we may
extend the result of Theorem \ref{cor1} to arbitrary joint moments.
For $k\geq1$, define
%
%
%e9 #&#
\begin{equation}
\label{WqTASEP} W^{k}_{\geq0}= \bigl\{\vec{n}
=(n_1,n_2,\ldots,n_k)\in(
\mathbb{Z}_{>0})^k\dvtx n_1\geq n_2
\geq \cdots\geq n_k\geq0 \bigr\}.
\end{equation}

\begin{figure}

\includegraphics{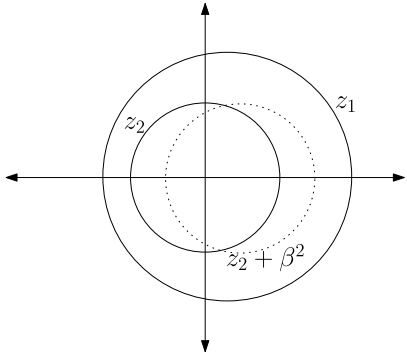}

\caption{Valid contours for equation (\protect\ref{1.6}) with $k=2$. The
inner contour is $z_2$ and the $z_1$ contour contains the image of
$z_2$ plus $\beta^2$.}\label{nested}
\end{figure}

%
%th1.7 #&#
\begin{theorem}\label{semiprop}
For all $k\geq1$ and $\vec{n}\in W^{k}_{\geq0}$,
%
%
%e10 #&#
\begin{equation}
\label{1.6} \mathbb{E} \Biggl[\prod_{i=1}^{k}
Z_\beta(t,n_i) \Biggr] = \frac
{1}{(2\pi
\iota)^k} \oint\cdots
\oint\prod_{1\leq a<b\leq k} \frac{z_a
-z_b}{z_a -z_b -\beta^2} \prod
_{i=1}^k \frac
{e^{t(z_i-1)}}{z_i^{n_i}}\frac{dz_i}{z_i},\hspace*{-33pt}
\end{equation}
where the integration contour for $z_a$ is a closed curve containing 0,
and the image under addition by $\beta^2$ of the integration contours
for $z_b$ for all $b>a$ (\,for an illustration of possible contours see
Figure~\ref{nested}).
\end{theorem}

%
%f1 #&#

This theorem is proved in Section~\ref{semidisc}. Asymptotics of this
formula yield all the moment Lyanpunov exponents. By Brownian scaling
it suffices to consider just $\beta=1$.

%
%f2 #&#
\begin{figure}[t]

\includegraphics{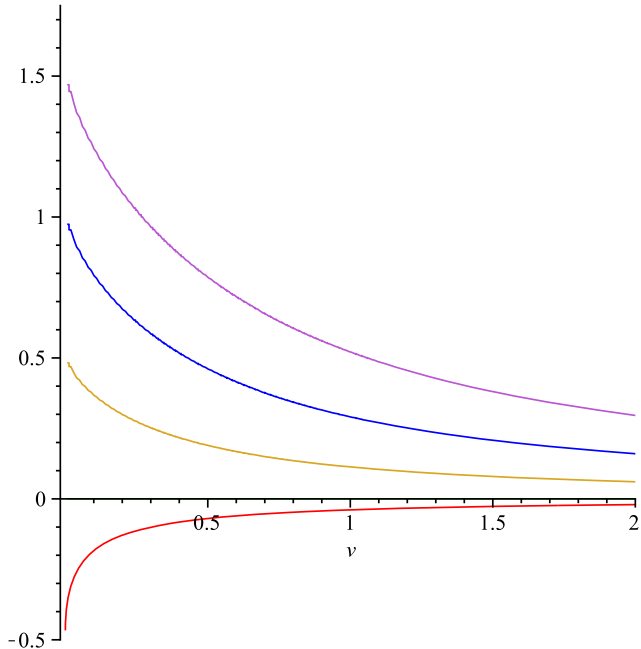}

\caption{Plot of one-sided parabolic Anderson model Lyapunov exponents
versus the velocity $\nu$ (which plays the role of the diffusion
constant). Here we have normalized $Z_1(t,\lfloor\nu t\rfloor)$ by
the zero noise solution, so that $\gamma_1(1;\nu)=0$. The lowest
curve in the plot is the normalized $\tilde\gamma_1(1;\nu)$, and the
higher curves are the normalized $\gamma_k(1;\nu) /k$ (increasing in
height with $k$). This demonstrates the intermittency of this parabolic
Anderson model, and the shape of this plot is very similar to those on
page 105 of \cite{CarMol}.}\label{lyapunov_plot_modified2}
\end{figure}

%
%th1.8 #&#
\begin{theorem}\label{kpam}
For any $k\geq1$ and $\nu>0$, the $k$th moment Lyapunov exponent at
velocity $\nu$ for the one-sided ($p=2$ and $q=0$) nearest-neighbor
parabolic Anderson model with $\beta=1$ is given by
\[
\gamma_k(1;\nu) = H_k \bigl(z^0_k
\bigr),
\]
where
\[
H_k(z) = \frac{k(k-3)}{2} + k z - \nu\log \Biggl(\prod
_{i=0}^{k-1} (z+i) \Biggr)
\]
and where $z^0_k$ is the unique solution to $H_k'(z)=0$ with $z\in
(0,\infty)$.
\end{theorem}
This theorem is proved in Section~\ref{semidisc} via asymptotics of
Theorem \ref{semiprop}. Figure~\ref{lyapunov_plot_modified2} records
the plot of the various Lyapunov exponents.

Note that in the one-sided case, the almost sure Lyapunov exponent
defined in~(\ref{tildegamma}) was conjectured in \cite{OCon-Yor} and
proved in \cite{OConnellMoriarty}, and it is given by (for $\beta=1$)
\[
\tilde\gamma_1(1;\nu) = -\frac{3}{2} + \inf
_{t>0} \bigl(t-\nu\Psi(t) \bigr),
\]
where $\Psi(t):=[\log\Gamma]'(t)$ is the digamma
function.\vadjust{\goodbreak}

There are two ideas which are behind the results of this paper. The
first idea is the content of Propositions \ref{uprop} and
\ref{prop:semiODEs} which show that one can compute the moments of the
parabolic Anderson model via solving a system of coupled ODEs with
spatial variables $\vec{n}\in W^{k}_{\geq0}$ and specific boundary
conditions. This reduction to solving ODEs on $W^{k}_{\geq0}$ only
works for $k=1,2$ with the general $p,q$ nearest-neighbor model.
However, for $p=2$ and $q=0$, the reduction holds for all $k$. The
second idea is that the system of ODEs can be explicitly solved via a
certain nested-contour integral ansatz that originated from \cite
{BorCor}. This is the content of Propositions~\ref{thm1}
and~\ref{solveonesided}.

%s1.4 #&#
\subsection{Outline}
The rest of the paper is as follows: in Section~\ref{nnpamsec} we show
how the moments of the parabolic Anderson model can be computed via a
coupled system of ODEs. We then solve this system and use this solution
to prove Theorems~\ref{cor1}~and~\ref{secondlyap}. In
Section~\ref{semidisc} we show how in the one-sided model, all moments
can be computed via ODEs, and we provide integral formulas which solve
these ODEs. From this we are able to prove Theorems \ref{semiprop} and
\ref{kpam}. In the \hyperref[app]{Appendix} we include a nonrigorous
replica trick calculation (used extensively in the physics literature)
and show how from this calculation one recovers the almost sure
Lyapunov exponent for the one-sided model; we also briefly discuss the
continuous space parabolic Anderson model (i.e., the stochastic heat
equation with multiplicative noise) and record its moments and Lyapunov
exponents.

%s2 #&#
\section{Nearest-neighbor parabolic Anderson model}\label{nnpamsec}

The first step in our computation of the moment Lyapunov exponents of
the parabolic Anderson model is the following reduction to a coupled
system of ODEs with two-body delta interaction. Recall the definition
of the nearest-neighbor parabolic Anderson model $Z_{\beta}(t,n)$ and
the operator $\Delta^{p,q}$ given in the \hyperref[sec1]{Introduction}.
Write $ [\Delta^{p,q} ]_i$ for the operator which acts as $\Delta
^{p,q}$ on the $i$th spatial coordinate.

%
%pr2.1 #&#
\begin{proposition}\label{prop1}
Assume $v\dvtx \mathbb{R}_{+}\times\mathbb{Z}^k\to\mathbb{R}$ solves:
\begin{longlist}[(3)]
\item[(1)] for all $\vec{n}\in\mathbb{Z}^k$ and $t\in\mathbb{R}_{+}$,
\[
\frac{d}{dt} v(t;\vec{n}) = \mathbf{H}v(t;\vec{n}), \qquad\mathbf{H}=
\frac
{1}{2}\sum_{i=1}^{k} \bigl[
\Delta^{p,q} \bigr]_i + \frac{1}{2} \beta^2
\mathop{\sum_{a,b=1}}_{a\neq b}^{k}
\mathbf{1}_{n_a=n_b};
\]

\item[(2)] for all permutations of indices $\sigma\in S_k$, $v(t;\sigma
    \vec{n}) = v(t;\vec{n})$;

\item[(3)] for all $\vec{n}\in\mathbb{Z}^k$, $\lim_{t\to0} v(t;\vec
    {n}) = \prod_{i=1}^{k} \mathbf{1}_{n_i=0}$;

\item[(4)] for all $T>0$, there exists $c,C>0$ such that for all $\vec
    {n}\in\mathbb{Z}^k$ and all $t\in[0,T]$,
\[
\bigl|v(t;\vec{n})\bigr| \leq c e^{C \|n\|_1}.
\]
\end{longlist}
Then for $\vec{n}\in\mathbb{Z}^k$, $v(t;\vec{n}) = \mathbb{E}
[\prod_{i=1}^{k} Z_\beta(t,n_i) ]$.
\end{proposition}

\begin{pf}
This result is well known and can be found, for instance, in
Proposition 6.1.3 of \cite{BorCor}. The purpose of the fourth
hypothesis on $v$ is to ensure uniqueness of solutions to the system of
ODEs given by the first three hypotheses. The fact that this
exponential growth hypothesis is sufficient for uniqueness can be
proved in the same manner as given in the proof of Proposition~4.9 in
\cite{BCS}.

One way to see why this should be true is to consider the Feynman--Kac
representation for $Z_{\beta}(t,n)$ which is given in equation
(\ref{FK}). The $k$ factors of $Z_{\beta}$ lead to~$k$ paths. The
expectation $\mathbb{E}$ over the Gaussian disorder (white-noise) can
be taken inside the path expectations $\mathcal{E}$ and calculated
exactly yielding the exponential of the pair-wise local time for the
$k$ paths. This accounts for the delta interaction seen above.
\end{pf}

It is a priori not clear how one would start to solve the system of
ODEs in the above proposition, one reason being that it is
inhomogeneous in space. An idea from integrable systems (related to the
coordinate Bethe Ansatz) is to instead try to solve a homogeneous
system of ODEs and put the inhomogeneity into a~boundary condition. If
the number of boundary conditions is $k-1$, then there is generally
hope in solving the system by combining fundamental solutions of the
homogeneous system in such a way that the initial data and boundary
conditions are met.

For the general $p,q$ case, it appears that this reduction to $k-1$
boundary conditions only works when $k=2$ (in which case there is just
one boundary condition). When $p=2$ and $q=0$ the reduction works for
all $k$; see Section~\ref{semidisc}.\looseness=1

%
%pr2.2 #&#
\begin{proposition}\label{uprop}
Assume $u\dvtx \mathbb{R}_{+}\times\mathbb{Z}^2\to\mathbb{R}$ solves:
\begin{longlist}[(3)]
\item[(1)] For all $\vec{n}\in\mathbb{Z}^2$ and $t\in\mathbb{R}_{+}$,
\[
\frac{d}{dt} u(t;\vec{n}) = \frac{1}{2}\sum
_{i=1}^{2} \bigl[\Delta^{p,q}
\bigr]_i u(t;\vec{n});
\]

\item[(2)] For $\vec{n}$ such that $n_1=n_2=n$
\begin{eqnarray*}
T_{\beta} u(t;\vec{n})&:=&\beta^2 u(t;n,n) + \frac{p}{2}
u(t;n,n-1) + \frac{q}{2} u(t;n+1,n)
\\
&&{}- \frac{p}{2} u(t;n-1,n) -
\frac{q}{2} u(t;n,n+1) = 0;
\end{eqnarray*}

\item[(3)] For all $\vec{n}\in\mathbb{Z}^k$ such that $n_1\geq n_2$,
    $\lim_{t\to 0} u(t;\vec{n}) = \prod_{i=1}^{2} \mathbf{1}_{n_i=0}$;

\item[(4)] For all $T>0$, there exists $c,C>0$ such that for all $\vec
    {n}\in\mathbb{Z}^k$ such that \mbox{$n_1\geq n_2$} and all $t\in[0,T]$,
\[
\bigl|u(t;\vec{n})\bigr| \leq c e^{C \|n\|_1}.
\]
\end{longlist}
Then for $\vec{n}\in\mathbb{Z}^2$ such that $n_1\geq n_2$, $u(t;\vec
{n}) =
v(t;\vec{n}) = \mathbb{E} [\prod_{i=1}^{2} Z_\beta
(t,n_i) ]$.
\end{proposition}

\begin{pf}
We show that restricted to $n_1\geq n_2$, $u(t;\vec{n})$ symmetrically
extended to $\mathbb{Z}^2$ solves the system of ODEs in Proposition
\ref{prop1} and hence $u(t;\vec{n}) = v(t;\vec{n})$. For $n_1>n_2$ it
is clear that $u$ and $v$ solve the same equation. For $n_1=n_2=n$,
\begin{eqnarray*}
&& \frac{d}{dt} u(t;n,n)
\\
&&\qquad = \frac{1}{2} \bigl( pu(t;n-1,n) + qu(t;n+1,n)
\\
&&\hspace*{43pt}{} + pu(t;n,n-1) + q u(t;n,n+1) - 4u(t;n,n) \bigr)
\\
&&\qquad = \bigl(\beta^2-2 \bigr) u(t;n,n) + pu(t;n,n-1) + q u(t;n+1,n),
\end{eqnarray*}
where the second line followed from the relation imposed by the
assumption (2). Now compare this to the equation $v(t;n,n)$ that satisfies:
\begin{eqnarray*}
&& \frac{d}{dt} v(t;n,n)
\\
&&\qquad = \frac{1}{2} \bigl( pv(t;n-1,n) + q
v(t;n+1,n) - 2 v(t;n,n) \bigr)
\\
&&\quad\qquad{} +\frac{1}{2} \bigl(pv(t;n,n-1) + q v(t;n,n+1) - 2 v(t;n,n) \bigr) +
\beta^2 v(t;n,n)
\\
&&\qquad = \bigl(\beta^2-2 \bigr) v(t;n,n) + p v(t;n,n-1) + qv(t;n+1,n),
\end{eqnarray*}
where the second line followed from the symmetry hypothesis on $v$.
Observe that on the diagonal $n_1=n_2$, both $u$ and $v$ solve the same
equation. Therefore they solve the same equation for all $n_1\geq n_2$
and hence (since the other hypotheses of Proposition \ref{prop1} are
satisfied) $u(t;\vec{n}) = v(t;\vec{n})$.
\end{pf}

We may now explicitly solve the system of ODEs defined in Proposition~\ref{uprop}.

%
%pr2.3 #&#
\begin{proposition}\label{thm1}
For $k=2$ and $n_1\geq n_2$, the system of ODEs in
Proposition~\ref{uprop} is uniquely solved by
%
%
%e11 #&#
\begin{eqnarray}\label{eqn1}
u(t;n_1,n_2) &=& \frac{1}{(2\pi\iota)^2} \oint
\oint\frac{(pz_1-qz_1^{-1}) - (pz_2 - q z_2^{-1})}{ (pz_1-qz_1^{-1}) - (pz_2 - q
z_2^{-1}) - 2\beta^2}
\nonumber\\[-8pt]\\[-8pt]
&&\hspace*{53pt}{}\times  F^{p,q}_{t,n_1}(z_1)
F^{p,q}_{t,n_2}(z_2) \frac
{dz_1}{z_1}\frac{dz_2}{z_2},\nonumber
\end{eqnarray}
where
\[
F^{p,q}_{t,n}(z) = z^{-n} e^{(t/2)(pz+qz^{-1} -2)}
\]
and where the contour of $z_1$ is the unit circle and the contour for
$z_2$ is a circle around 0 of radius sufficiently small so as not to
include any poles of the integrand aside from $z_2=0$; see Figure~\ref{nested}.
\end{proposition}

%
%re2.4 #&#
\begin{remark}
The right-hand side of (\ref{eqn1}) is easy to generalize to all $k$;
see, for example, Section~6.1.2 of \cite{BorCor} or Proposition
\ref{solveonesided} below in the totally asymmetric case where $p=2$
and $q=0$. However, it is not at all clear if such a generalization
would have anything to do with $\mathbb{E} [\prod_{i=1}^{k}
Z_\beta(t,n_i) ]$.
\end{remark}

Before proving this proposition, we note that Theorem \ref{cor1}
follows as an immediate corollary of the above result and Proposition
\ref{uprop}.

%By virtue of Proposition \ref{uprop}, the solution for $u(t;n_1,n_2)$
%given in (\ref{eqn1}) is equal to the desired two-point joint moment.

%For $n_1\geq n_2$,
%(pz_1-qz_1^{-1}) - (pz_2 - %q z_2^{-1}) - 2\beta^2} F_{t,n_1}(z_1)
%F_{t,n_2}(z_2) \frac{dz_1}{z_1} \frac{dz_2}{z_2},
%with $F_{t,n}(z)$ and the contours of integration specified in Theorem

\begin{pf*}{Proof of Proposition \ref{thm1}}
We prove this proposition by checking the hypotheses of Proposition
\ref{uprop}. Hypothesis~1 follows from the Leibniz rule and the
observation that
\[
\frac{d}{dt} F^{p,q}_{t,n}(z) = \frac{1}{2}
\Delta^{p,q} F^{p,q}_{t,n}(z).
\]

To check hypothesis 2, we apply $T_\beta$ to $u(t;\vec{n})$ when
$n_1=n_2=n$. The operator $T_\beta$ can be taken inside the
integration. It acts on $ F^{p,q}_{t,n}(z_1) F^{p,q}_{t,n}(z_2)$ as
\begin{eqnarray*}
&& T_\beta \bigl(F^{p,q}_{t,n}(z_1) F^{p,q}_{t,n}(z_2) \bigr)
\\
&&\qquad =
-\frac {1}{2} F^{p,q}_{t,n}(z_1) F^{p,q}_{t,n}(z_2) \bigl(
\bigl(pz_1-qz_1^{-1} \bigr) - \bigl(pz_2-qz_2^{-1} \bigr) -2 \beta^2
\bigr).
\end{eqnarray*}
The factor $ ((pz_1-qz_1^{-1}) - (pz_2-qz_2^{-1}) -2\beta^2 )$
cancels with the same term in the denominator of the integrand, yielding
\begin{eqnarray*}
T_\beta u(t;\vec{n}) &=& \frac{1}{(2\pi\iota)^2}
\oint\oint \bigl( \bigl(pz_1-qz_1^{-1} \bigr) - \bigl(pz_2 - q z_2^{-1}
\bigr) \bigr)
\\
&&\hspace*{52pt}{}\times F^{p,q}_{t,n}(z_1) F^{p,q}_{t,n}(z_2)
\frac{dz_1}{z_1} \frac{dz_2}{z_2}.
\end{eqnarray*}
Since the integration contours are the same for both $z_1$ and $z_2$
and since $ ((pz_1-qz_1^{-1}) - (pz_2 - q z_2^{-1}) )$ is
antisymmetric in $z_1$ and $z_2$, while the rest of the integrand is
symmetric, one immediately sees that the integral is zero as desired to
check hypothesis 2.

Hypothesis 3 is checked via residue calculus. The $t\to0$ limit can be
taken inside the integrand, and we are left to show that
for $n_1\geq n_2$,
%
%
%e12 #&#
\begin{equation}
\label{eqnabove} \frac{1}{(2\pi\iota)^2} \oint\oint C(z_1,z_2)
z_1^{-n_1}z_2^{-n_2} \frac{dz_1}{z_1}
\frac{dz_2}{z_2} = \mathbf{1}_{n_1=0}\mathbf{1}_{n_2=0},
\end{equation}
where
\[
C(z_1,z_2) = \frac{ (pz_1-qz_1^{-1}) - (pz_2 - q z_2^{-1})}{
(pz_1-qz_1^{-1}) - (pz_2 - q z_2^{-1}) - 2\beta^2}
\]
and where the contour of $z_1$ is the unit circle and the contour for
$z_2$ is a circle around 0 of radius sufficiently small so as not to
include any poles of the integrand aside from $z_2=0$.
\begin{longlist}[(3)]
\item[(1)] If $n_2<0$, then in (\ref{eqnabove}) we may shrink the $z_2$
    contour to 0. Observe that for $z_1$ fixed on the specified
    contour, $C(z_1,z_2)$ is analytic in $z_2$ in a small neighborhood
    of $z_2=0$, with a value of $C(z_1,0)=1$. The term $z_2^{-n_2}
    \frac{dz_2}{z_2}$ does not have a pole at~0 (because \mbox{$n_2<0$}) and
    hence in this case $u(0;\vec{n})=0$.

\item[(2)] If $n_2=0$, then in (\ref{eqnabove}) let us shrink the $z_2$
    contour to 0. The term $z_2^{-n_2} \frac{dz_2}{z_2}$ has a simple
    pole at 0 and hence the integral evaluates as
\[
\frac{1}{2\pi\iota} \oint z_1^{-n_1}\frac{dz_1}{z_1}=
\mathbf{1}_{n_1=0}.
\]

\item[(3)] If $n_2>0$ this implies that $n_1>0$ as well. Then we can
    expand $z_1$ to infinity. As we do this, we encounter a pole at
    $z_1$ such that
\[
\bigl(pz_1-q z_1^{-1} \bigr) -
\bigl(pz_2-qz_2^{-1} \bigr) - 2
\beta^2=0.
\]
For each $z_2$ there is only one such pole $r(z_2)$ which comes
when $z_1\approx\frac{q}{p} z_2^{-1}$ for small $|z_2|$. The reason
why only one pole is crossed is because the other pole coming from
this term is approximately $z_2$, which is already contained inside
the $z_1$ contour. Before analyzing the residue, observe that
because $n_1>0$, there is at least quadratic decay in $z_1$ at
infinity, so there is no pole at infinity. Thus, the integral in
$z_1$ is given by its negative residue at $r(z_2)$.

The negative residue at $z_1=r(z_2)$ is evaluated as
%
%
%e13 #&#
\begin{equation}
\label{entint} -\frac{1}{2\pi\iota} \oint\frac{2\beta^2 r(z_2)}{pr(z_2) +
q(r(z_2))^{-1}} \bigl(r(z_2)
\bigr)^{-n_1-1} z_2^{-n_2-1} \,dz_2,
\end{equation}
where the integral in $z_2$ is over a small circle around the origin.
It is easy to see that
\[
\frac{2\beta^2 r(z_2)}{pr(z_2) + q(r(z_2))^{-1}}
\]
is analytic in a neighborhood of $z_2=0$ and its value at $z_2=0$ is
$2\beta^2/p$. Thus the entire integral (\ref{entint}) can be
evaluated as the residue at $z_2=0$. Since $n_1\geq n_2$, there is, in
fact, no pole at $z_2=0$, thus the integral equals 0.
\end{longlist}
Combining the above cases we see that the only case in which $u(0;\vec
{n})$ is nonzero is when $n_1=n_2=0$, in which case it is 1. This
confirms the initial data of hypothesis~3.

Hypothesis 4 follows via easy bounds of the integrand of
$u(t;\vec{n})$.
\end{pf*}

Having proved the two-point moment formula for the nearest-neighbor
para\-bolic Anderson model, we can now extract the second moment
Lyapunov exponent via asymptotic analysis.

\begin{pf*}{Proof of Theorem \ref{secondlyap}}
We will present a complete proof only in the case $p=q=1$ since this
simplifies the (rather technical) analysis. For the moment we keep the
$p$ and $q$ and only set them equal when necessary.

Let us start by proving $\gamma_1=0$ from the formula
\[
\mathbb{E} \bigl[Z_\beta(t,n) \bigr] = \frac{1}{2\pi\iota} \oint
F^{p,q}_{t,n}(z) \frac{dz}{z},
\]
which one easily checks via either Proposition \ref{prop1} or
\ref{uprop}. Let $n=\lfloor(p-q)t\rfloor$ and observe that [up to an
insignificant correction coming from the fractional difference between
$n$ and $(p-q)t$]
\[
\mathbb{E} \bigl[Z_\beta(t,n) \bigr] = \frac{1}{2\pi\iota} \oint
e^{t
G(z)} \frac{dz}{z}, \qquad G(z) = pz+qz^{-1} -2 -
(p-q)\log z.
\]
We want to study this as $t\to\infty$; thus we can use the standard
Laplace method (see Lemma \ref{asymptoticslemma} for $\ell=1$) to
perform the asymptotics. The critical point equation for $G(z)$ is
\[
G'(z) = p - qz^{-2} -(p-q)z^{-1} = 0,
\]
which is solved by $z=1$ or $z=-q/p$. The critical point $z=1$
corresponds to the larger value of $G(z)$, namely $G(1)=0$. Observe
that we can deform the $z$ contour to lie on the unit circle $e^{\iota
\theta}$. As a function of $\theta$ along this contour
$\operatorname{Re} [G(z(\theta))] = 2 \cos(\theta) - 2$. This shows
that the real part of $G(z)$ decays monotonically with respect to the
angle $\theta$ away from $z=1$. In the vicinity of $z=1$,
$\operatorname{Re}[G(z)]$ decays quadratically in the imaginary
directions. Invoking Lemma \ref{asymptoticslemma} for $\ell=1$ shows
that
\[
\gamma_1:=\lim_{t\to\infty} \frac{1}{t} \log
\bigl(\mathbb{E} \bigl[Z_\beta \bigl(t,(p-q)t \bigr) \bigr] \bigr) =0.
\]

To calculate $\gamma_2$ we use the formula for $\mathbb{E}
[Z_{\beta
}(t,n_1)Z_{\beta}(t,n_2) ]$ proved in Theorem \ref{cor1}. In
order to perform the asymptotics we would like to deform the $z_2$
contour to coincide with the $z_1$ contour. While doing this we
encounter a pole and must take into account the associated residue, in
addition to the evaluation of the remaining integral on the new
contours. The result of this manipulation is
\[
\mathbb{E} \bigl[Z_{\beta} \bigl(t,(p-q)t \bigr)^2 \bigr] =
\mathrm{(A)} + \mathrm{(B)},
\]
where
\begin{eqnarray*}
\mathrm{(A)} &=& \frac{1}{(2\pi\iota)^2} \oint\oint\frac{
(pz_1-qz_1^{-1}) - (pz_2 - q z_2^{-1})}{ (pz_1-qz_1^{-1}) - (pz_2 - q
z_2^{-1}) - 2\beta ^2}
\\
&&\hspace*{53pt}{}\times F^{p,q}_{t,(p-q)t}(z_1) F^{p,q}_{t,(p-q)t}(z_2)
\frac{dz_1}{z_1} \frac{dz_2}{z_2}
\end{eqnarray*}
with the $z_2$ contour coinciding with the $z_1$ contour,
and
%
%
%e14 #&#
\begin{equation}
\label{Beqn} \mathrm{(B)} = \frac{1}{2\pi\iota} \oint\frac{2\beta^2}{p+q
[s(z_2)]^{-2}} e^{tH(z_2)}
\frac{dz_2}{z_2}.
\end{equation}
Note that this residue term comes from $z_1=s(z_2)$ where $s(z_2)$ is
given in the statement of the theorem. %The residue at $z_1=s(z_2)$ is
%then evaluated as above.
Since the definition of $s(z)$ involves a square-root, for $z$ complex
we specify that for $z=re^{\iota\theta}$ with
$\theta\in(-\infty,\infty)$, $\sqrt{z} = \sqrt{r} e^{\iota\theta/2}$.

As it will turn out, the residue term \textup{(B)} has a larger
exponential growth rate than the integral term \textup{(A)} and hence
accounts entirely for the value of $\gamma_2$. To see this, we compute
the exponent for both terms. We claim that
\[
\lim_{t\to\infty} \frac{1}{t} \log \bigl[\mathrm{(A)} \bigr] = 0.
\]
It is easy to see why this is true. The contours for $z_1$ and $z_2$
can be chosen so that there exists a positive constant $C$ such that
along the contours
\[
\bigl\llvert g(z_1,z_2) \bigr\rrvert< C,
\]
where
\[
g(z_1,z_2)= \frac{ (pz_1-qz_1^{-1}) - (pz_2 - q z_2^{-1})}{
(pz_1-qz_1^{-1}) - (pz_2 - q z_2^{-1}) - 2\beta^2}.
\]
Likewise, in a neighborhood of $(z_1,z_2)=(1,1)$, one checks that
\[
c\bigl|z_1^2 z_2 - z_2 -
z_1 z_2^2 +z_1\bigr| \leq \bigl
\llvert g(z_1,z_2) \bigr\rrvert
\]
for a small, yet positive $c$. One easily checks the remaining
assumptions necessary to apply Lemma \ref{asymptoticslemma} (with
$\ell=2$) and therefore finds that as $t\to\infty$, the growth of the
integral defining \textup{(A)} is governed by the value of
$G(z_1)+G(z_2)$ at the critical point $(1,1)$. By comparison to the
calculation performed above for $\gamma_1$, we find that
\[
\lim_{t\to\infty} \frac{1}{t} \log \bigl[\mathrm{(A)} \bigr] = 2
\gamma_1.
\]
Since $\gamma_1=0$, the claimed result for \textup{(A)}
follows.\vadjust{\goodbreak}

Turn now to the residue term \textup{(B)} and call $z_2$ simply $z$.
From this point on we will assume that $p=q=1$ to simplify the
analysis. A similar, albeit lengthy, analysis can be performed for all
$p$ and $q$. Notice that
\[
z^0_2 = \frac{1}{2} \bigl(-\beta^2+
\sqrt{4+\beta^4} \bigr)
\]
is a critical point of $H_2(z)$ and that the contour for $z$ can be
deformed without crossing any singularities of (\ref{Beqn}) to the
contour $\Gamma$ parameterized by $z= z^0_2 e^{\iota\theta}$ for
$\theta\in[0,2\pi]$. We wish to use Laplace's method by applying Lemma
\ref{asymptoticslemma} with $\ell=1$. It is straightforward to check
hypotheses~(1), (3) and (4) of the lemma. Hypothesis~(2) requires more
work.

To check hypothesis (2a) of Lemma \ref{asymptoticslemma} observe that
$H_2'(z^0_2)=0$, $H_2''(z^0_2)\neq0$ and $H_2(z)$ is analytic in a
neighborhood of $z^0_2$. Thus it immediately follows that it behaves
locally like $H(z^0_2)+c (z-z^0_2)^2 + o((z-z^0_2)^{2})$. Hypothesis
(2b) requires that
\[
\rho(\theta):= 2 \operatorname{Re} \bigl[H_2(z)-H_2
\bigl(z^0_2 \bigr) \bigr]
\]
(the factor of $2$ is irrelevant) is strictly negative for all $z\in
\Gamma\setminus\{z^0_2\}$. In fact, by symmetry of $H$ through the
real axis, $\rho(\theta) = \rho(2\pi-\theta)$ and hence this
strict negativity needs only be checked for $\theta\in(0,\pi]$. By
utilizing the fact that
\begin{eqnarray*}
\bigl(z^0_2 \bigr)^{-1} &=& \frac{1}{2}
\bigl(\beta^2+\sqrt{4+\beta^4} \bigr),
\\
s(z)^{-1} &=&
\frac{-(z-z^{-1}+2\beta^2) + \sqrt{(z-z^{-1}+2\beta ^2)^2 + 4}}{2}
\end{eqnarray*}
we find that
\begin{eqnarray*}
\rho(\theta) &=& \sqrt{4+\beta^4} \bigl(\cos(\theta)-2 \bigr)
\\
&&{} + \operatorname{Re} \bigl[\bigl(\beta^4 \bigl(2-\cos(\theta) \bigr)^2-
\bigl(4+\beta^4 \bigr)\sin(\theta)^2
\\
&&\hspace*{30pt}{}  + 4 + \iota2 \beta^2
\bigl(2-\cos(\theta) \bigr)\sqrt{4+ \beta^4} \sin( \theta)\bigr)^{1/2} \bigr].
\end{eqnarray*}
Since for $a,b$ real,
\[
\operatorname{Re}\sqrt{a+\iota b} = \sqrt{ \frac{a+\sqrt{a^2+b^2}}{2}},
\]
checking hypothesis (2b) reduces to checking that for $\theta\in
(0,\pi]$
\[
\rho(\theta) = \sqrt{4+\beta^4} \bigl(\cos(\theta)-2 \bigr) + \sqrt{
\frac
{a+\sqrt{a^2+b^2}}{2}}<0,
\]
where
\begin{eqnarray*}
a &=&\beta^4 \bigl(2-\cos(\theta) \bigr)^2- \bigl(4+ \beta^4
\bigr)\sin(\theta)^2 + 4,
\\
b&=& 2 \beta^2 \bigl(2-\cos(\theta)
\bigr)\sqrt{4+\beta^4} \sin( \theta).
\end{eqnarray*}

It is straight-forward to check that $\rho(\pi)<0$ hence by
continuity of $\rho(\theta)$ it suffices to prove that $\rho(\theta
)\neq0$ for $\theta\in(0,\pi]$. A simple calculation shows that
this is equivalent to showing that
\[
128 \bigl(4+\beta^4 \bigr) \bigl(\cos(\theta)-2 \bigr)^2
\sin(\theta/2)^2 \neq0
\]
on $\theta\in(0,\pi]$ which is immediately verified. Thus we have
proved hypothesis (2b) of Lemma \ref{asymptoticslemma}.

%This inequality is equivalent (via elementary algebraic manipulations)
%to showing that
%$$
%-16 (4 + \beta^4) \Big(-36 - \beta^4 + 32 \cos(\theta) + (-4 +
%$$
%Letting $x=\cos(\theta)$ the above inequality is equivalent to showing
%that
%f(x)&=&-2 (1 - x^2) \beta^4 (4 + \beta^4) - \Big(4 + (-2 + x)^2
%- (1 - x^2) (4 + \beta^4)\Big)^2>0
%for $x\in[-1,1)$. The above equation is of fifth degree in $x$,
%however one checks that $x=1$ is a root. Direct inspection shows that
%all four other roots are off the real axis (in fact they can be
%computed explicitly since after factoring $(1-x)$ the equation is
%quintic). Thus, to show the desired inequality, it suffices to show
%that the derivative of the right-hand side of (\ref{fineqn}) at $x=1$
%is negative. Computing, one sees that $f'(1) = -4 (64 + 20 \beta^4+
%inequality, it immediately follows that $\rr(\theta)<0$ for all $\theta

Applying Lemma \ref{asymptoticslemma} to (\ref{Beqn}) we find that
\[
\lim_{t\to\infty} \frac{1}{t} \log\mathrm{(A)} = H_2
\bigl(z^0_2 \bigr).
\]
Since $H_2(z^0_2)$ is positive [as compared to the contribution from
the \textup{(A)} term asymptotics] we conclude that
$\gamma_2(\beta;0)=H_2(z^0_2)$.
\end{pf*}

%s3 #&#
\section{The one-sided parabolic Anderson model}\label{semidisc}

We now focus on the one-sided case where $p=2$ and $q=0$. Recall the
definition of $W^{k}_{\geq0}$ given in (\ref{WqTASEP}) and the
definition of $\Delta^{p,q}$ given in (\ref{deltadef}). In
particular, $\Delta^{2,0}f(n) = f(n-1)-f(n)$.% is the one-sided
%derivative operator.

%
%pr3.1 #&#
\begin{proposition}\label{prop:semiODEs}
Fix $k\in\mathbb{Z}_{>0}$.
\begin{longlist}[\textup{Part} 1.]
\item[\textup{Part} 1.] Assume $v\dvtx
    \mathbb{R}_{+}\times(\mathbb{Z}_{\geq-1})^k$ solves:
\begin{longlist}[(3)]
\item[(1)] for all $\vec{n}\in(\mathbb{Z}_{\geq-1})^k$ and
    $t\in\mathbb {R}_{+}$,
\[
\frac{d}{dt} v(t;\vec{n}) = \mathbf{H}v(\tau;\vec{n}), \qquad\mathbf{H}=
\frac{1}{2}\sum_{i=1}^{k} \bigl[
\Delta^{2,0} \bigr]_i + \frac
{1}{2}
\beta^2 \mathop{\sum_{a,b=1}}_{a\neq b}^k
\mathbf{1}_{n_a=n_b};
\]

\item[(2)] for all permutations of indices $\sigma\in S_k$,
    $v(t;\sigma \vec{n}) = v(t;\vec{n})$;

\item[(3)] for all $\vec{n}\in W^{k}_{\geq0}$, $v(0;\vec{n}) =
    \prod_{i=1}^{k} \mathbf{1}_{n_i=0}$;

\item[(4)] for all $\vec{n}\in(\mathbb{Z}_{\geq-1})^k$ such that
    $n_k=-1$, $v(t;\vec{n}) \equiv0$ for all $t\in\mathbb{R}_{+}$.
\end{longlist}
Then for all $\vec{n}\in W^{k}_{\geq0}$, $\mathbb{E} [\prod_{i=1}^{k}
Z_{\beta}(t,n_i) ] = v(t;\vec{n})$.

\item[\textup{Part} 2.] Assume $u\dvtx
    \mathbb{R}_{+}\times(\mathbb{Z}_{\geq-1})^k
    \to\mathbb{R}$ solves:
\begin{longlist}[(3)]
\item[(1)] for all $\vec{n}\in(\mathbb{Z}_{\geq-1})^k$ and
    $t\in\mathbb {R}_{+}$,
\[
\frac{d}{dt } u(t;\vec{n}) = \frac{1}{2}\sum
_{i=1}^{k} \bigl[\Delta^{2,0}
\bigr]_i u (t;\vec{n});
\]

\item[(2)] for all $\vec{n}\in(\mathbb{Z}_{\geq-1})^k$ such that
    for some $i\in \{1,\ldots, k-1\}$, $n_i=n_{i+1}$,
\[
\bigl( \bigl[\Delta^{2,0} \bigr]_i - \bigl[
\Delta^{2,0} \bigr]_{i+1} - 2\beta^2 \bigr) u(t;
\vec{n}) = 0;
\]

\item[(3)] for all $\vec{n}\in(\mathbb{Z}_{\geq-1})^k$ such that
    $n_k=-1, u(t;\vec{n}) \equiv0$ for all $t \in\mathbb{R}_{+}$;

\item[(4)] for all $\vec{n}\in W^{k}_{\geq0}$, $u(0;\vec{n}) =
    \prod_{i=1}^{k} \mathbf{1}_{n_i=0}$.
\end{longlist}
Then for all $\vec{n}\in W^{k}_{\geq0}$, $\mathbb{E} [\prod_{i=1}^{k}
Z_{\beta}(t,n_i) ] = u(t;\vec{n})$.
\end{longlist}
\end{proposition}

\begin{pf}
This is contained in Proposition 6.3 of \cite{BCS}.
\end{pf}

%
%re3.2 #&#
\begin{remark}
Part~1 of the proposition is essentially a specialization of
Proposition \ref{prop1} to the case $p=2$, $q=0$. However, due to the
delta initial data and the one-sided nature of the operator $\Delta
^{2,0}$, $v(t;\vec{n})$ with $\vec{n}\in\{-1,0,\ldots, m\}^k$ evolves
autonomously as a closed system of ODEs. This ensures uniqueness of
solutions and explains why we no longer require the at-most exponential
growth hypothesis which was present in Proposition \ref{prop1}. Part~2
of the proposition is an extension of Proposition \ref{uprop} to
general $k$, but only for $p=2$, $q=0$. The fact that this holds for
all $k$ is what enables us to solve for higher than second moments in
this one-sided case.
\end{remark}

The system of ODEs in part~2 of Proposition \ref{prop:semiODEs} can be
solved via a nested-contour integral ansatz introduced in \cite
{BorCor} and further developed in \cite{BCS}. This yields the following
generalization of Proposition \ref{thm1} to all $k$ but only for the
one-side ($p=2$, $q=0$) case.

%
%pr3.3 #&#
\begin{proposition}\label{solveonesided}
For all $k\geq1$ and $\vec{n}\in W^{k}_{\geq0}$, the system of ODEs in
part~2 of Proposition \ref{prop:semiODEs} is uniquely solved by
%
%
%e15 #&#
\begin{equation}
\label{eqn12} u(t;\vec{n}) = \frac{1}{(2\pi\iota)^k} \oint\cdots\oint \prod
_{1\leq a<b\leq k} \frac{z_a -z_b}{z_a -z_b -\beta^2} \prod_{i=1}^k
\frac{e^{t(z_i-1)}}{z_i^{n_i}}\frac{dz_i}{z_i},
\end{equation}
where the integration contour for $z_a$ is a closed curve containing 0
and the image under addition by $\beta^2$ of the integration contours
for $z_b$ for all $b>a$.
\end{proposition}

Before proving this proposition, we note that Theorem \ref{semiprop}
follows as an immediate corollary of the above result and Proposition
\ref{prop:semiODEs}, part~2.

\begin{pf*}{Proof of Proposition \ref{solveonesided}}
We check the four conditions for the system of ODEs in part~2 of
Proposition \ref{prop:semiODEs}. Condition~(1) follows by Leibnitz rule
and the fact that
\[
\frac{d}{dt} \frac{e^{t(z-1)}}{z^{n}} = \frac{1}{2}\Delta^{2,0}
\frac{e^{t(z-1)}}{z^{n}}
\]
for all $z\in\mathbb{C}\setminus\{0\}$.\vadjust{\goodbreak}

Condition~(2) follows by applying $ ( [\Delta^{2,0} ]_i - [\Delta^{2,0}
]_{i+1} - 2\beta^2 )$ to the integrand of the right-hand side of
(\ref{eqn12}). The effect of this operator is to bring out an extra
factor of $2(z_i-z_{i+1}-\beta^2)$. This factor cancels the
corresponding term in the denominator of the product over $a<b$.
Without the pole associated with this term, it is possible to deform
the $z_i$ and $z_{i+1}$ contours to coincide, and since $n_i=n_{i+1}$,
we find that
\[
\bigl( \bigl[\Delta^{2,0} \bigr]_i - \bigl[
\Delta^{2,0} \bigr]_{i+1} - 2\beta^2 \bigr) u(t;
\vec{n}) = \oint dz_{i}\oint \,dz_{i+1} (z_i-z_{i+1})
f(z_i)f(z_{i+1}),
\]
where $f(z)$ includes all of the other integrations aside from those in
$z_i$ and $z_{i+1}$. The above integral is clearly 0 by skew-symmetry,
thus confirming condition~(2) as desired.

Condition~(3) follows by observing that if $n_k=-1$, then on the
right-hand side of (\ref{eqn12}) there is no pole at 0 in the $z_k$
variable. By Cauchy's theorem, this means that since the $z_k$ contour
only contains 0 and no other poles of the integrand, the entire
integral is 0, as desired.

Condition~(4) follows from three easy residue calculations. From above,
if \mbox{$n_k<0$}, the integral in (\ref{eqn12}) is zero. Similarly, if
$n_1>0$ the integrand in (\ref{eqn12}) has no pole at infinity and the
$z_1$ contour can be freely deformed to infinity. By Cauchy's theorem
this means that the entire integral is zero. The only possible nonzero
value of $u(0;\vec{n})$ is (due to the ordering of the elements in
$\vec{n}\in W^{k}_{\geq0}$) when\vadjust{\goodbreak} $n_1=\cdots=
n_k=0$. The value of $u$ for this choice of $\vec{n}$ is readily
calculated via residues to equal 1, just as desired.
\end{pf*}

We now show how asymptotics of the result of Theorem \ref{semiprop}
yield a proof of the moment Lyapunov exponents claimed in Theorem \ref{kpam}.

\begin{pf*}{Proof of Theorem \ref{kpam}}
The starting point of this proof is the moment formula of Theorem
\ref{semiprop}. Setting all $n_i\equiv\lfloor\nu t\rfloor$ and $\beta
=1$ we find that
%
%
%e16 #&#
\begin{eqnarray}\label{wefindthat}
&& \mathbb{E} \Biggl[\prod_{i=1}^{k}
Z_\beta \bigl(t,\lfloor\nu t \rfloor \bigr) \Biggr]
\nonumber\\[-11pt]\\[-11pt]
&&\qquad  = \frac{1}{(2\pi\iota)^k} \oint\cdots\oint\prod_{1\leq a<b\leq
k} \frac{z_a -z_b}{z_a -z_b -1} \prod_{i=1}^k
F^{2,0}_{t,\lfloor\nu
t\rfloor}(z_i) \frac{dz_i}{z_i},\nonumber
\end{eqnarray}
where the integration contour for $z_a$ is a closed curve containing 0
and the image under addition by $1$ of the integration contours for
$z_b$ for all $b>a$. From now on we will study the right-hand side of
the\vspace*{-1pt} above equality, with $F^{2,0}_{t,\lfloor\nu
t\rfloor}(z_i)$ replaced by $F^{2,0}_{t, \nu t}(z_i)$, as the
asymptotic effect of this modification is easily seen to be
inconsequential.

%The first step in our analysis is to turn the integration contours
%into infinite contours. We claim that equation (\ref{wefindthat})
%holds with the integral contour for $z_a$ given by $\{\alpha_a + x e^{
%where $\alpha_a> \alpha_b+1$ for all $b>a$. This claim is easily shown
%due to the fact that $|F_{t,n}^{2,0}(z)|=e^{t(\Re(z)-1)}$. This bound
%and Cauchy's theorem enables us to open the original closed contours
%into the desired infinite contours.

In order to perform the asymptotic analysis necessary to compute the
moment Lyapunov exponents, we would like to deform our contours to all
coincide so as to apply Lemma \ref{asymptoticslemma}. This requires
deforming all contours to pass through a specific critical point of
$\log F^{2,0}_{t, \nu t}(z)$. However, due to the nesting of the
contours, such a~deformation requires passing a number of poles. The
following lemma records the effect of such a deformation.

%
%f3 #&#
\begin{figure}%[b]

\includegraphics{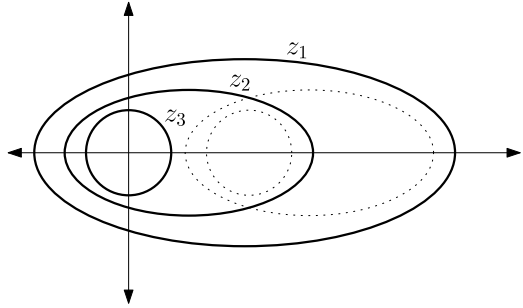}

\caption{Valid contours for Lemma \protect\ref{reslemma} with $k=3$. The
inner contour is $z_3$ and contains~0; the next contour is $z_2$ and
contains 0 and the image of the $z_3$ contour plus one; the outer
contour is $z_1$ and contains 0 and the image of the $z_2$ and $z_3$
contours plus one. The images of the $z_3$ and $z_2$ contours after
adding one are indicated by the dotted lines.}\label{plus1nested}
\end{figure}

%
%le3.4 #&#
\begin{lemma}\label{reslemma}
Consider a function $f(z)$ which is analytic in $\mathbb{C}\setminus\{
0\}$. For $k\geq1$, set\vspace*{-1pt}
\[
\mu_k= \frac{1}{(2\pi\iota)^k} \oint\cdots\oint\prod
_{1\leq
A<B\leq k} \frac{z_A-z_B}{z_A-z_B-1} \frac{f(z_1)\cdots
f(z_k)}{z_1\cdots z_k} \,dz_1
\cdots dz_k,
\]
where the integration contour for $z_A$ contains 0 and the image under
addition by~$1$ of the integration contours for $z_B$ for all $B>A$;
see Figure~\ref{plus1nested}. Then
\[
\mu_k= k! \sum_{\lambda=1^{m_1}2^{m_2}\cdots\vdash k} I_{\lambda},
\]
where
%
%
%e17 #&#
\begin{eqnarray}\label{Ilambda}
I_\lambda &=& \frac{1}{m_1!m_2!\cdots(2\pi\iota)^{\ell
(\lambda)}}\nonumber
\\[-2pt]
&&{}\times \oint\cdots\oint\det \biggl[ \frac{1}{\lambda_i+w_i-w_j}
\biggr]_{i,j=1}^{\ell(\lambda)}
\\[-2pt]
&&\hspace*{46pt}{}\times \prod_{j=1}^{\ell(\lambda)} f(w_j)f(w_j+1)
\cdots f(w_j+\lambda_j-1) \,dw_j.\nonumber
\end{eqnarray}
Here $\lambda=1^{m_1}2^{m_2}\cdots\vdash k$ means $\lambda$ is a
partition $\lambda= (\lambda_1\geq\lambda_2\geq\cdots\geq0)$ such that
$\sum\lambda_i =k$; $m_i$ records the number of entries of $\lambda$
equal to $i$; $\ell(\lambda)$ is the number of nonzero entries in
$\lambda$; and for $1\leq j\leq\ell(\lambda)$ the $w_j$ contours are
all chosen to be the same contour as $z_k$.
\end{lemma}

\begin{pf}
This is given in \cite{BorCor}, Proposition 6.2.7, and is proved by
taking a scaling limit (with $q\mapsto e^{-\varepsilon}$, $z=\mapsto
e^{-\varepsilon z}$, $w\mapsto e^{-\varepsilon w}$ and
$\varepsilon\to0$) of \cite{BorCor}, Proposition~3.2.1. Alternatively
another proof is given in \cite{BorCorhalfspace} as Proposition 5.1.
\end{pf}

The following lemma will also be helpful in completing our asymptotics.

%
%le3.5 #&#
\begin{lemma}\label{36}
Consider $I_{\lambda}$ in (\ref{Ilambda}) with $f(z) = F^{2,0}_{t,\nu
t}(z)$ and write $I_{\lambda}(t)$ to emphasize the $t$ dependence.
Then
\[
\gamma_\lambda= \lim_{t\to\infty} \frac{1}{t} \log
I_\lambda(t)
\]
exists and is given by\vspace*{-1pt}
%
%
%e18 #&#
\begin{equation}
\label{glam} \gamma_\lambda= \sum_{j=1}^{\ell(\lambda)}
\gamma_{\lambda_j},
\end{equation}
where
\[
\gamma_{r} = H_{r} \bigl(z^0_r
\bigr).
\]
Here, as in the statement of Theorem \ref{kpam},
\[
H_{r}(z) = \frac{r(r-3)}{2} + r z - \nu\log \Biggl(\prod
_{i=0}^{r-1} (z+i) \Biggr)
\]
and $z^0_r$ is the unique solution to $H_{r}'(z)=0$ with $z\in
(0,\infty)$.
\end{lemma}

\begin{pf}
First observe that we can take the contours of integration for $w_j$
in~$I_{\lambda}$ to be a large circle containing $\{0,-1,\ldots,
-\lambda_j+1\}$. This is because before having applied the identity in
Lemma \ref{reslemma}, we could take the $z_j$ contours to be large
enough nested circles so that $z_k$ contains $0,-1,\ldots,
-\lambda_j+1$.

Next we observe that the integrals defining $I_{\lambda}$ match the
form of (\ref{It}) in Lemma~\ref{asymptoticslemma} with $\ell=\ell
(\lambda)$ and (using $w$'s instead of $z$'s)
\[
g(w_1,\ldots, w_\ell) = \frac{1}{m_1!m_2!\cdots} \det \biggl[
\frac
{1}{\lambda_i+w_i-w_j} \biggr]_{i,j=1}^{\ell(\lambda)}
\]
and $G_{j}(w_j) = H_{\lambda_j}(w_j)$. If we can show that the four
hypotheses of Lem\-ma~\ref{asymptoticslemma} apply, then the result
claimed in the present lemma follows immediately.\looseness=1

By convexity, $G_j(w_j)$ has exactly one critical point along $w_j\in
(0,\infty)$. Call this point $w^0_j$. Without changing the value of the
integrals, we can freely deform the contour of integration for $w_j$ to
a contour $\Gamma_j$ which is defined (see Figure~\ref{semicircle} for
an illustration) as the union of a long vertical line segment going
through $w^0_j$ and a semi-circle enclosing $\{ 0,-1,\ldots,
-\lambda_j+1\}$, with radius large enough so that for all
$r\in\{0,-1,\ldots, -\lambda_j+1\}$ and all $w_j\in\Gamma
_j\setminus\{w^0_j\}$, $|w_j-r|>|w^0_j-r|$. For this choice of contour
it is clear that all of the hypotheses of Lemma \ref{asymptoticslemma}
hold. In particular hypothesis (2a) holds since
$\operatorname{Re}(w_j)$ is constant along the vertical portion of the
contour and decreasing on the circular part; and
$\operatorname{Re}[\log(w_j+i)]>\operatorname{Re}[\log(w^0_j+i)]$ for
all $w_j\neq w^0_j$ along $\Gamma_j$.
\end{pf}

We can now complete the proof of Theorem \ref{kpam}. From
straightforward comparison of growth of exponentials it follows that
\[
\gamma_k:= \lim_{t\to\infty} \frac{1}{t} \log
\biggl(\sum_{\lambda\vdash k} I_\lambda(t) \biggr) = \max
_{\lambda\vdash k} \gamma_\lambda.
\]
Observe that by combining equation (\ref{wefindthat}) with Lemmas
\ref{reslemma} and \ref{36} we find that
%
%
%e19 #&#
\begin{equation}
\label{gmax} \gamma_k(1;\nu)= \max_{\lambda\vdash k}
\gamma_{\lambda},
\end{equation}
where $\gamma_k(1;\nu)$ is defined in (\ref{gamma}), and $\gamma
_{\lambda}$ is defined above in (\ref{glam}). We claim now that this
maximum is attained for all $k$ at $\gamma_{k}$ and hence $\gamma
_k(1;\nu)=\gamma_k$. If we can show this, then the theorem is proved.

%
%f4 #&#
\begin{figure}

\includegraphics{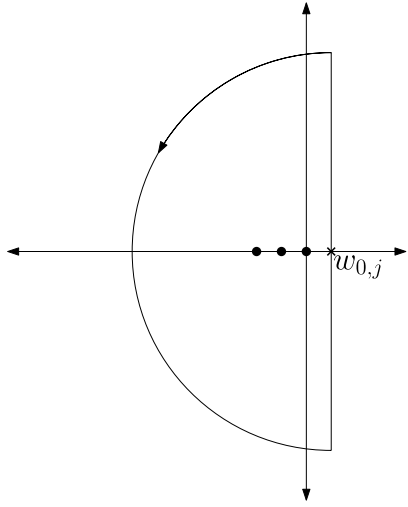}

\caption{The contour $\Gamma_j$ is a vertical line segment with real
part $w_{0,j}$ joined with a semi-circle which encloses $\{0,-1,\ldots,
-\lambda_j+1\}$ as well has the property that $r\in\{0,-1,\ldots,
-\lambda_j+1\}$ and all $w_j\in\Gamma_j\setminus\{w_{0,j}\}$,
$|w_j-r|>|w_{0,j}-r|$. In words, this means that the distance from the
points on~$\Gamma_j$ to the various elements of the set $\{0,-1,\ldots,-\lambda_j+1\}$ is minimized for $w_{0,j}$ and strictly larger
otherwise.}\label{semicircle}
\end{figure}

For $k=1$, this result is immediate since there is only one partition
of $1$. For $k=2$ there are two partitions to consider $\lambda=(1,1)$
and $\lambda=(2)$. We claim that for all $\nu>0$, $2\gamma_1-\gamma
_2>0$. This can be proved by explicit evaluation. Observe that
\begin{eqnarray*}
\gamma_1 &=& -1 + \nu+\nu\log\nu,
\\
\gamma_2 &=& -2 + \nu+
\sqrt{1+\nu^2} - \nu\log \biggl(\frac{1}{2} \nu \bigl(\nu+
\sqrt{1+\nu^2} \bigr) \biggr).
\end{eqnarray*}
The difference $f(\nu):= \gamma_2-2\gamma_1$ goes to zero as $\nu$
goes to infinity and has derivative $\log(2\nu) - \log(\nu+ \sqrt{1+\nu
^2})$ which is negative for all $\nu>0$. This shows that $f(\nu
)>0$ for all $\nu>0$ and hence
\[
\gamma_2(1;\nu)= \max_{\lambda\vdash2} \gamma_{\lambda}
= \gamma_2.
\]

Note that we have now shown that $\gamma_1(1;\nu)< \gamma_2(1;\nu)/2$
and hence we can apply Lemma \ref{inter} to show the intermittency of
all of the moment Lyapunov exponents. We now proceed by induction on
$k$. Assume that for all $j\leq k$, we have proved that
$\gamma_k(1;\nu)= \gamma_k$. As a base case we have $k=1$ and $2$. By
intermittency we know that for any partition of $k+1$ aside from
$\lambda=(k+1)$, $\gamma_{k+1}(1;\nu)$ must strictly exceed $\sum_{i}
\gamma_{\lambda_i}(1;\nu)$ which, by induction, equals $\sum_{i}
\gamma_{\lambda_i}$ as well. By (\ref{gmax}) this implies that the
maximum over $\lambda\vdash k+1$ must be attained for the partition
$\lambda= (k+1)$ and hence $\gamma_{k+1}(1;\nu) = \gamma_{k+1}$ as
desired to prove the inductive step and complete the proof of
Theorem~\ref{kpam}.
\end{pf*}

%sA #&#
\begin{appendix}\label{app}
\section*{Appendix}

In the first subsection of this Appendix we show how (using our moment
Lyapunov exponents for the one-side parabolic Anderson model) the
physics replica trick leads nonrigorously to the correct formula for
the almost sure Lyapunov exponent $\tilde\gamma_1$. This almost sure
exponent is already-known rigorously \cite{OConnellMoriarty}, so the
below calculation should be thought of as a nontrivial check of the
efficacy of the replica trick.

In the second subsection of this Appendix we apply the nested contour
integral methods to compute all of the moment Lyapunov exponents for
the continuum parabolic Anderson model (i.e., stochastic heat equation
with multiplicative noise). These exponents have been known for some
time and were first computed in the physics literature by Kardar \cite
{K} and then in the math literature by Bertini and Cancrini's; see
\cite{BC} Theorem 2.6 and remark after it.

The final subsection of this Appendix contains a version of Laplace's
method for computing asymptotics of integrals (and was referenced
earlier in the paper).

%sA.1 #&#
\subsection{The replica trick}\label{replica}
The replica trick is an idea which goes back to Kac~\cite{Kac} and
which has received a great deal of attention within the statistical
physics community. In its most basic form, one hopes to extract the
almost sure Lyapunov exponent from the knowledge of all of the moment
Lyapunov exponents. The reader should be warned that what follows is
extremely nonrigorous. However, we include it since it is a validation
of the replica trick in the context of the one-sided parabolic Anderson
model; see \cite{K} for this approach implemented in the continuous
model discussed below in Section~\ref{SHE}.

We would like to compute the almost sure Lyapunov exponent
\[
\tilde\gamma_1(1;\nu):= \lim_{t\to\infty}
\frac{1}{t} \mathbb{E} \bigl[\log Z_{1}(t,\nu t) \bigr].
\]
Note that even though we have taken the expectation of $\log Z_1(t,\nu
t)$, this should not affect the value of the almost sure exponent.
Recall that for $z\in\mathbb{C}\setminus\mathbb{R}_{-}$,
%
%
%eA.1 #&#
\begin{equation}
\label{logz} \log z = \lim_{k\to0} \frac{z^k-1}{k}.
\end{equation}
We have shown in Theorem \ref{kpam} that
\[
\mathbb{E} \bigl[Z_1(t,\nu t)^k \bigr] \approx
e^{t \gamma_k(1;\nu)} = e^{t H_k(z^0_k)}
\]
for
\[
H_k(z) = \frac{k(k-3)}{2} + kz - \nu\log \Biggl(\prod
_{i=0}^{k-1} (z+i) \Biggr) = \frac{k(k-3)}{2} + kz -
\nu\log\frac{ \Gamma
(z+k)}{\Gamma(z)}
\]
and $z^0_k$ is the unique minimum of $H_k(z)$ for $z\in(0,\infty)$.
This second expression has a clear analytic extension in $k$.

By using (\ref{logz}) and interchanging the two limits (without
justification) we have
\[
\tilde\gamma_1(1;\nu) = \lim_{t\to\infty}
\frac{1}{t} \lim_{k\to0} \frac{e^{t H_k(z^0_k)} -1}{k}.
\]
Notice that for $k$ near 0, $e^{t H_k(z^0_k)}\approx1 + t H_k(z^0_k)$, hence
\[
\tilde\gamma_1(1;\nu) = \lim_{t\to\infty}
\frac{1}{t} \lim_{k\to0} \frac{ t H_k(z^0_k)}{k}.
\]
The limit in $t$ can now be taken (since the factors of $t$ cancel),
and the limit in $k$ is achieved via L'H\^{o}pital's rule,
\[
\lim_{k\to0} \frac{H_k(z)}{k} = -\frac{3}{2} +z - \nu
\Psi(z),
\]
where $\Psi(z) = [\log\Gamma]'(z)$ is the digamma function. This
limit should be evaluated at the unique infimum over $z\in(0,\infty)$
and hence
\[
\tilde\gamma_1(1;\nu) = -\frac{3}{2} + \inf
_{z>0} \bigl(z-\nu\Psi(z) \bigr).
\]

This nonrigorous calculation does yield the proved value; cf.
\cite{OCon-Yor} and \cite{OConnellMoriarty}.

%sA.2 #&#
\subsection{The space--time continuum parabolic Anderson model}\label{SHE}

Consider the solution $\mathcal{Z}_{\beta}\dvtx \mathbb{R}_{+}\times
\mathbb{R}\to \mathbb{R}_{+}$ to the multiplicative stochastic heat
equation with delta function initial data,
%
%
%eA.2 #&#
\begin{equation}
\frac{d}{dt} \mathcal{Z}_{\beta} = \frac{1}{2} \Delta\mathcal
{Z}_{\beta} + \beta\dot{W} \mathcal{Z}_{\beta}, \qquad\mathcal
{Z}_{\beta}(0,x)=\delta_{x=0},
\end{equation}
where $\Delta$ is the Laplacian on $\mathbb{R}$, and $\delta_{x=0}$ is
the Dirac delta function. The solution $\mathcal{Z}_{\beta}(t,x)$ can
be thought of as the partition function for a space--time continuous
directed polymer in a white-noise environment \cite{ACQ}. In fact,
under a particular scaling, the parabolic Anderson model considered in
the previous sections, converges to the SHE \cite{QRMF}. The following
formula for joint moments can be found by applying this limit
transition to Theorem \ref{semiprop}. For all $k\geq1$ and all $x_1\leq
x_2\leq\cdots\leq x_k$ in~$\mathbb{R}^k$,
%
%
%eA.3 #&#
\begin{eqnarray}\label{SHEprop}
&& \mathbb{E} \Biggl[\prod_{i=1}^{k}
\mathcal{Z}_{\beta}(t,x_i) \Biggr]
\nonumber\\[-8pt]\\[-8pt]
&&\qquad = \frac{1}{(2\pi\iota)^k} \oint
\cdots\oint\prod_{1\leq a<b\leq k} \frac{z_a -z_b}{z_a -z_b -\beta^2} \prod
_{i=1}^k e^{(t/2)z_i^2 + x_i z_i}
\,dz_i,\nonumber
\end{eqnarray}
where $z_j\in\beta^2 \alpha_j + \iota\mathbb{R}$ and $\alpha_1
>\alpha
_2+1>\alpha_3+2>\cdots> \alpha_k + k-1$.

The formula of Theorem \ref{semiprop} solved the system of ODEs in
Proposition \ref{prop:semiODEs}. Likewise, (\ref{SHEprop}) solves a PDE
which is called the ``quantum delta Bose gas''; see Section~6.2 of
\cite{BorCor} for a discussion on this, as well as remarks on certain
gaps in a~rigorous statement to this effect.

From (\ref{SHEprop}) it is possible to compute the moment Lyapunov
exponents for the SHE (we restrict attention now to $x_i\equiv0$ and
to $\beta=1$ since general $x_i\equiv x$ and $\beta$ can be achieved
from the resulting formula via Brownian scaling). The result is that
\[
\gamma_k:= \lim_{t\to\infty} \frac{1}{t} \log
\bigl(\mathbb{E} \bigl[\mathcal{Z}(t,0)^k \bigr] \bigr) =
\frac{k^3-k}{24}.
\]
This reproduces Kardar's formula \cite{K} and agrees with Bertini and
Cancrini's Theorem 2.6 and remark after it; see \cite{BC}.

This calculation is done by deforming all contours to the imaginary
axis and considering the growth in $t$ of the various residue terms. As
in Theorem \ref{kpam}, the Lyapunov exponent comes from the
(ground-state) term when $z_1=z_2+1 = z_3+2 = \cdots=z_k +k-1$. There
remains only one free variable of integration in this residue term, and
the main part of the integrand is the following exponential:
\begin{eqnarray*}
&& \exp \biggl\{\frac{t}{2} \bigl(z^2 + (z+1)^2 +
\cdots+ (z+k-1)^2 \bigr) \biggr\}
\\
&&\qquad =
\exp \biggl\{\frac{kt}{2} \biggl(z+\frac{k-1}{2}
\biggr)^2 + \frac{t}{2} \biggl(1+ 2^2 + \cdots
(k-1)^2 - k\frac{(k-1)^2}{4} \biggr) \biggr\}.
\end{eqnarray*}
Deforming the $z$-contour to $\iota\mathbb{R}-\frac{k-1}{2}$ shows that
this term behaves like
\[
\exp \biggl\{\frac{t}{2} \biggl(1+ 2^2 +
\cdots(k-1)^2 - k\frac
{(k-1)^2}{4} \biggr) \biggr\} = \exp \biggl\{t
\biggl(\frac
{k^3-k}{24} \biggr) \biggr\}
\]
from which the result readily follows.

%sA.3 #&#
\subsection{Laplace's method}

The following lemma is a version of Laplace's method for computing
asymptotics of integrals. The proof is an easy modification of the
usual proof of Laplace's method \cite{Erd}.

\setcounter{theorem}{0}
%
%leA.1 #&#
\begin{lemma}\label{asymptoticslemma}
Consider a contour integral
%
%
%eA.4 #&#
\begin{equation}
\label{It} I(t) = \frac{1}{(2\pi\iota)^\ell} \oint_{\Gamma_1} \cdots \oint
_{\Gamma_\ell} g(z_1,\ldots,z_\ell) \exp \Biggl(t \sum
_{j=1}^{\ell
} G_j(z_j)
\Biggr) \,dz_1 \cdots dz_\ell.
\end{equation}

Assume that:
\begin{longlist}[(3)]
\item[(1)] For each $j$, $\Gamma_j$ is a closed piecewise smooth
    contour.

\item[(2)] For each $j$, there exists $z^0_j\in\Gamma_j$ such that:
\begin{enumerate}[a]
\item[(a)] $\operatorname{Re}[G_j(z)]
    <\operatorname{Re}[G_j(z^0_j)]$ for all $z\in \Gamma_j$ not
    equal to $z^0_j$;

\item[(b)] $G_j'(z^0_j)=0$ and in a neighborhood of $z^0_j$,
    $G_j(z) = G_j(z^0_j) + c(z-z^0_j)^r + o((z-z^0_j)^{r})$ for
    some $r\geq2$.
\end{enumerate}

\item[(3)] There exists a (nonidentically zero) rational function
    $R(z_1,\ldots, z_\ell)$ such that in a neighborhood of
    $(z^0_1,\ldots, z^0_j)$,
\[
\bigl|R(z_1,\ldots, z_\ell)\bigr| \leq\bigl|g(z_1,\ldots,
z_\ell)\bigr|.
\]

\item[(4)] There exists a positive constant $C$ such that for all
    $z_j\in \Gamma_j$ ($j=1,\ldots, \ell$),
\[
\bigl|g(z_1,\ldots, z_\ell)\bigr| \leq C.
\]
\end{longlist}

Then
\[
\lim_{t\to\infty} \frac{1}{t} \log I(t) = \sum
_{j=1}^{\ell} G_j \bigl(z^0_j
\bigr).
\]
\end{lemma}
\end{appendix}

\section*{Acknowledgments} We thank Frank den Hollander, Davar
Koshnevisan and Quentin Berger for discussions on the parabolic
Anderson model and useful comments on this article.


%suskaldyti doi

% imsref loaded by linak, 2013-11-20 12:54:24

\printaddresses


\begin{thebibliography}{22}
% BibTex style file: ims.bst, 2013-01-28
% Default style options (sort=0,type=number).
% Used options (sort=1,type=number).
%b1 #&#
\bibitem{ACQ}
\begin{barticle}[mr]
\bauthor{\bsnm{Amir},~\bfnm{Gideon}\binits{G.}},
  \bauthor{\bsnm{Corwin},~\bfnm{Ivan}\binits{I.}} \AND
  \bauthor{\bsnm{Quastel},~\bfnm{Jeremy}\binits{J.}}
(\byear{2011}).
\btitle{Probability distribution of the free energy of the continuum directed
  random polymer in {$1+1$} dimensions}.
\bjournal{Comm. Pure Appl. Math.}
\bvolume{64}
\bpages{466--537}.
\bid{doi={10.1002/cpa.20347}, issn={0010-3640}, mr={2796514}}
\bptok{imsref}%
\end{barticle}
\endbibitem

%b2 #&#
\bibitem{Berger}
\begin{barticle}[mr]
\bauthor{\bsnm{Berger},~\bfnm{Quentin}\binits{Q.}} \AND
  \bauthor{\bsnm{Lacoin},~\bfnm{Hubert}\binits{H.}}
(\byear{2011}).
\btitle{The effect of disorder on the free-energy for the random walk pinning
  model: Smoothing of the phase transition and low temperature asymptotics}.
\bjournal{J. Stat. Phys.}
\bvolume{142}
\bpages{322--341}.
\bid{doi={10.1007/s10955-010-0110-x}, issn={0022-4715}, mr={2764128}}
\bptok{imsref}%
\end{barticle}
\endbibitem

%b3 #&#
\bibitem{BC}
\begin{barticle}[mr]
\bauthor{\bsnm{Bertini},~\bfnm{Lorenzo}\binits{L.}} \AND
  \bauthor{\bsnm{Cancrini},~\bfnm{Nicoletta}\binits{N.}}
(\byear{1995}).
\btitle{The stochastic heat equation: {F}eynman--{K}ac formula and
  intermittence}.
\bjournal{J. Stat. Phys.}
\bvolume{78}
\bpages{1377--1401}.
\bid{doi={10.1007/BF02180136}, issn={0022-4715}, mr={1316109}}
\bptok{imsref}%
\end{barticle}
\endbibitem

%b4 #&#
\bibitem{BorCor}
\begin{bmisc}[auto:STB|2013/10/14|10:36:11]
\bauthor{\bsnm{Borodin},~\bfnm{A.}\binits{A.}} \AND
  \bauthor{\bsnm{Corwin},~\bfnm{I.}\binits{I.}}
(\byear{2013}).
\bhowpublished{Macdonald processes. \textit{Probab. Theory Related Fields}. To
  appear}.
\bptok{imsref}%
\end{bmisc}
\endbibitem


%b5 #&#
\bibitem{BorCorhalfspace}
\begin{bmisc}[auto:STB|2013/10/14|10:36:11]
\bauthor{\bsnm{Borodin},~\bfnm{A.}\binits{A.}} \AND
  \bauthor{\bsnm{Corwin},~\bfnm{I.}\binits{I.}}
(\byear{2013}).
\bhowpublished{Directed random polymers via nested contour integrals. Unpublished manuscript}.
\bptok{imsref}%
\end{bmisc}
\endbibitem

%b6 #&#
\bibitem{BCF}
\begin{bmisc}[auto:STB|2013/10/14|10:36:11]
\bauthor{\bsnm{Borodin},~\bfnm{A.}\binits{A.}},
  \bauthor{\bsnm{Corwin},~\bfnm{I.}\binits{I.}} \AND
  \bauthor{\bsnm{Ferrari},~\bfnm{P.~L.}\binits{P.~L.}}
(\byear{2013}).
\bhowpublished{Free energy fluctuations for directed polymers in random media
  in $1+1$ dimension. \textit{Comm. Pure Appl. Math.} To appear.}
\bptok{imsref}%
\end{bmisc}
\endbibitem

%b7 #&#
\bibitem{BCS}
\begin{bmisc}[auto:STB|2013/10/14|10:36:11]
\bauthor{\bsnm{Borodin},~\bfnm{A.}\binits{A.}},
  \bauthor{\bsnm{Corwin},~\bfnm{I.}\binits{I.}} \AND
  \bauthor{\bsnm{Sasamoto},~\bfnm{T.}\binits{T.}}
(\byear{2013}).
\bhowpublished{From duality to determinants for \mbox{$q$-}TASEP and ASEP. \textit{Ann. Probab.} To appear}.
\bptok{imsref}%
\end{bmisc}
\endbibitem

%b8 #&#
\bibitem{CarMol}
\begin{barticle}[mr]
\bauthor{\bsnm{Carmona},~\bfnm{Ren{\'e}~A.}\binits{R.~A.}} \AND
  \bauthor{\bsnm{Molchanov},~\bfnm{S.~A.}\binits{S.~A.}}
(\byear{1994}).
\btitle{Parabolic {A}nderson problem and intermittency}.
\bjournal{Mem. Amer. Math. Soc.}
\bvolume{108}
\bpages{viii+125}.
\bid{doi={10.1090/memo/0518}, issn={0065-9266}, mr={1185878}}
\bptok{imsref}%
\end{barticle}
\endbibitem

%b9 #&#
\bibitem{CC}
\begin{barticle}[mr]
\bauthor{\bsnm{Comets},~\bfnm{Francis}\binits{F.}} \AND
  \bauthor{\bsnm{Cranston},~\bfnm{Michael}\binits{M.}}
(\byear{2013}).
\btitle{Overlaps and pathwise localization in the {A}nderson polymer model}.
\bjournal{Stochastic Process. Appl.}
\bvolume{123}
\bpages{2446--2471}.
\bid{doi={10.1016/j.spa.2013.02.010}, issn={0304-4149}, mr={3038513}}
\bptnote{check year}%
\bptok{imsref}%
\end{barticle}
\endbibitem

%b10 #&#
\bibitem{CSY}
\begin{bincollection}[mr]
\bauthor{\bsnm{Comets},~\bfnm{Francis}\binits{F.}},
  \bauthor{\bsnm{Shiga},~\bfnm{Tokuzo}\binits{T.}} \AND
  \bauthor{\bsnm{Yoshida},~\bfnm{Nobuo}\binits{N.}}
(\byear{2004}).
\btitle{Probabilistic analysis of directed polymers in a random environment: A
  review}.
In \bbooktitle{Stochastic Analysis on Large Scale Interacting Systems}
\bpages{115--142}.
\bpublisher{Math. Soc. Japan}, \blocation{Tokyo}.
\bid{mr={2073332}}
\bptok{imsref}%
\end{bincollection}
\endbibitem

%b11 #&#
\bibitem{ICReview}
\begin{barticle}[mr]
\bauthor{\bsnm{Corwin},~\bfnm{Ivan}\binits{I.}}
(\byear{2012}).
\btitle{The {K}ardar--{P}arisi--{Z}hang equation and universality class}.
\bjournal{Random Matrices: Theory Appl.}
\bvolume{1}
\bpages{1130001, 76}.
\bid{doi={10.1142/S2010326311300014}, issn={2010-3263}, mr={2930377}}
\bptok{imsref}%
\end{barticle}
\endbibitem

%b12 #&#
\bibitem{Fest}
\begin{bbook}[auto:STB|2013/10/14|10:36:11]
\beditor{\bsnm{Deuschel},~\bfnm{J.-D.}\binits{J.-D.},
\bsnm{Gentz},~\bfnm{B.}\binits{B.},
\bsnm{K\"onig},~\bfnm{W.}\binits{W.},
\bparticle{van}\bsnm{Renesse},~\bfnm{M.-K.}\binits{M.-K.},
\bsnm{Scheutzow},~\bfnm{M.}\binits{M.}
\textup{and}
\bsnm{Schmock},~\bfnm{U.}\binits{U.}}, eds.
(\byear{2012}).
\btitle{Probability in Complex Physical Systems. In Honour of Erwin Bolthausen and J\"urgen G\"artner}.
\bseries{Springer Proceedings in Mathematics}
\bvolume{11}.
\bpublisher{Springer}, \blocation{Berlin}.
\bptok{imsref}%
\end{bbook}
\endbibitem

%b13 #&#
\bibitem{Erd}
\begin{bbook}[mr]
\bauthor{\bsnm{Erd{\'e}lyi},~\bfnm{A.}\binits{A.}}
(\byear{1956}).
\btitle{Asymptotic Expansions}.
\bpublisher{Dover}, \blocation{New York}.
\bid{mr={0078494}}
\bptok{imsref}%
\end{bbook}
\endbibitem

%b14 #&#
\bibitem{denHoll}
\begin{barticle}[mr]
\bauthor{\bsnm{Greven},~\bfnm{A.}\binits{A.}} \AND \bauthor{\bparticle{den}
  \bsnm{Hollander},~\bfnm{F.}\binits{F.}}
(\byear{2007}).
\btitle{Phase transitions for the long-time behavior of interacting
  diffusions}.
\bjournal{Ann. Probab.}
\bvolume{35}
\bpages{1250--1306}.
\bid{doi={10.1214/009117906000001060}, issn={0091-1798}, mr={2330971}}
\bptok{imsref}%
\end{barticle}
\endbibitem

%b15 #&#
\bibitem{Kac}
\begin{bmisc}[auto:STB|2013/10/14|10:36:11]
\bauthor{\bsnm{Kac},~\bfnm{M.}\binits{M.}}
(\byear{1968}).
\bhowpublished{On certain Toeplitz-like matrices and their relation to the
  problem of lattice vibrations. \textit{Arkiv for det Fysiske seminar i
  Trondheim} \textbf{11}.}
\bptok{imsref}%
\end{bmisc}
\endbibitem

%b16 #&#
\bibitem{K}
\begin{barticle}[mr]
\bauthor{\bsnm{Kardar},~\bfnm{Mehran}\binits{M.}}
(\byear{1987}).
\btitle{Replica {B}ethe ansatz studies of two-dimensional interfaces with
  quenched random impurities}.
\bjournal{Nuclear Phys. B}
\bvolume{290}
\bpages{582--602}.
\bid{doi={10.1016/0550-3213(87)90203-3}, issn={0550-3213}, mr={0922846}}
\bptok{imsref}%
\end{barticle}
\endbibitem

%b17 #&#
\bibitem{QRMF}
\begin{bmisc}[auto:STB|2013/10/14|10:36:11]
\bauthor{\bsnm{Moreno~Flores},~\bfnm{G.}\binits{G.}},
  \bauthor{\bsnm{Remenik},~\bfnm{D.}\binits{D.}} \AND
  \bauthor{\bsnm{Quastel},~\bfnm{J.}\binits{J.}}
\bhowpublished{Unpublished manuscript}.
\bptok{imsref}%
\end{bmisc}
\endbibitem

%b18 #&#
\bibitem{OConnellMoriarty}
\begin{barticle}[mr]
\bauthor{\bsnm{Moriarty},~\bfnm{J.}\binits{J.}} \AND
  \bauthor{\bsnm{O'Connell},~\bfnm{N.}\binits{N.}}
(\byear{2007}).
\btitle{On the free energy of a directed polymer in a {B}rownian environment}.
\bjournal{Markov Process. Related Fields}
\bvolume{13}
\bpages{251--266}.
\bid{issn={1024-2953}, mr={2343849}}
\bptok{imsref}%
\end{barticle}
\endbibitem

%b19 #&#
\bibitem{OCon}
\begin{barticle}[mr]
\bauthor{\bsnm{O'Connell},~\bfnm{Neil}\binits{N.}}
(\byear{2012}).
\btitle{Directed polymers and the quantum {T}oda lattice}.
\bjournal{Ann. Probab.}
\bvolume{40}
\bpages{437--458}.
\bid{doi={10.1214/10-AOP632}, issn={0091-1798}, mr={2952082}}
\bptok{imsref}%
\end{barticle}
\endbibitem

%b20 #&#
\bibitem{OCon-Yor}
\begin{barticle}[auto:STB|2013/10/14|10:36:11]
\bauthor{\bsnm{O'Connell},~\bfnm{N.}\binits{N.}} \AND
  \bauthor{\bsnm{Yor},~\bfnm{M.}\binits{M.}}
(\byear{2001}).
\btitle{Brownian analogues of Burke's theorem}.
\bjournal{Stochastic Process. Appl.}
\bvolume{96}
\bpages{285--304}.
\bptok{imsref}%
\end{barticle}
\endbibitem

%b21 #&#
\bibitem{SeppValko}
\begin{barticle}[mr]
\bauthor{\bsnm{Sepp{\"a}l{\"a}inen},~\bfnm{Timo}\binits{T.}} \AND
  \bauthor{\bsnm{Valk{\'o}},~\bfnm{Benedek}\binits{B.}}
(\byear{2010}).
\btitle{Bounds for scaling exponents for a {$1+1$} dimensional directed polymer
  in a {B}rownian environment}.
\bjournal{ALEA Lat. Am. J. Probab. Math. Stat.}
\bvolume{7}
\bpages{451--476}.
\bid{issn={1980-0436}, mr={2741194}}
\bptok{imsref}%
\end{barticle}
\endbibitem
\end{thebibliography}
\end{document}